\documentclass[10pt,reqno]{amsart}
\usepackage{amsfonts,amsmath,amssymb,amsthm}
\usepackage{latexsym}
\usepackage{eucal}
\usepackage[dvips]{graphics}
\usepackage[latin1]{inputenc}

\numberwithin{equation}{section}

\newtheorem{Theorem}{Theorem}[section]
\newtheorem{Lemma}{Lemma}[section]
\newtheorem{Proposition}{Proposition}[section]
\newtheorem{Corollary}{Corollary}[section]
\newtheorem{Remark}{Remark}[section]

\def\F{\mathcal{F}}

\begin{document}
\title[Some Systems of Coupled KdV Equations]
{On The Local Well-Posedness for Some Systems of Coupled
KdV Equations}

\setlength{\baselineskip}{1.3\baselineskip}

\author[B. Alvarez-Samaniego  
        and 
        X. Carvajal]{}

\email{balvarez@math.u-bordeaux1.fr, balvarez@math.uic.edu}
\email{carvajal@ime.unicamp.br}

%\subjclass[2000]{35A07, 35Q53.}

%\keywords{Hirota-Satsuma system, Gear-Grimshaw system, KdV equation.}
\thanks{$^1$Initially supported by FAPESP/Brazil 
grant No. 2002/02522-0, and later by CNRS/France.}
\thanks{$^2$Supported by FAPESP under grant No. 2002/08920-7.}
\thanks{{\it Date}: October 2005.}

\maketitle

\centerline{\scshape   Borys Alvarez-Samaniego$^{^1}$}
 \medskip

  {\footnotesize \centerline{Universit\'e Bordeaux I; IMB and CNRS UMR 5251}
  \centerline{351 Cours de la Lib\'eration, 33405 Talence Cedex, France}} 
 \medskip

\centerline{\scshape  Xavier Carvajal$^{^2}$}
 \medskip

 {\footnotesize \centerline{ Department of Mathematics, IMECC-UNICAMP}
 \centerline{C.P. 6065, CEP 13083-970, Campinas, SP, Brazil}}

\begin{abstract}
Using the theory developed by Kenig, Ponce, and Vega, we prove that
the Hirota-Satsuma system is locally well-posed in Sobolev spaces 
$H^s(\mathbb{R}) \times H^{s}(\mathbb{R})$ for $3/4<s\le1$.  
We introduce some Bourgain-type spaces $X_{s,b}^a$ for 
$a\not =0$, $s,b \in \mathbb{R}$  to obtain local well-posedness
for the Gear-Grimshaw system in $H^s(\mathbb{R})\times H^s(\mathbb{R})$
for $s>-3/4$, by establishing new mixed-bilinear estimates involving
the two Bourgain-type spaces  $X_{s,b}^{-\alpha_-}$ and
$X_{s,b}^{-\alpha_+}$ adapted to $\partial_t+\alpha_-\partial_x^3$ and
$\partial_t+\alpha_+\partial_x^3$ respectively, where 
$|\alpha_+|=|\alpha_-|\not = 0$. 
\hfill \break \vspace{ 4 mm}\break {\it {MSC:}} 35A07; 35Q53
\hfill \break \vspace{ 4 mm}\break {\it {Keywords:}} Hirota-Satsuma 
system; Gear-Grimshaw system; KdV equation 
\end{abstract}

%%%%%%%%%%%%%%%%%%%%%%%%%%%%%%%%%%%%%%%%%%%%%%%%%%
%  INTRODUCTION
%%%%%%%%%%%%%%%%%%%%%%%%%%%%%%%%%%%%%%%%%%%%%%%%%%

\section{Introduction}
In this paper we are concerned with two systems of coupled
KdV equations, namely the Hirota-Satsuma system and the
Gear-Grimshaw system.

First we consider local well-posedness (LWP) and ill-posedness of
the initial value problem (IVP) for the following system:
\begin{equation}\label{eq:hs}
  \left \{
   \begin{array}{l}
     u_t-a(u_{xxx}+6uu_x)=2bvv_x, \\
     v_t+v_{xxx}+3uv_x=0, \\
     u(0)=u_0, \;\; v(0)=v_0,
   \end{array}
  \right.
\end{equation}
known as the Hirota-Satsuma system which was
introduced in \cite{hs:hs} to describe the interaction of two long
waves with different dispersion relations.  Here $a, b$ are real
constants, and $u, v$ are real-valued functions of the two
real variables $x$ and $t$.  System (\ref{eq:hs}) is a set
of coupled Korteweg-de Vries (abbreviated KdV henceforth) equations,
and it is a generalization of the KdV equation (which is obtained
when $v=0$).  The Cauchy problem associated to (\ref{eq:hs}),
for the real and periodic case, was previously studied by P. F. He \cite{p:p},
for $b>0$, $-1<a<0$, and considering Sobolev indices $s \ge 3$.
It deserves remark that system (\ref{eq:hs}) has the following
conserved quantities:
\begin{eqnarray}
   V(u,v) &=& \int_{-\infty}^{+\infty} \Big{(} \frac{1+a}{2} u_x^2
   +bv_x^2-(1+a)u^3-buv^2\Big{)}dx, \label{eq:cons1} \\
   F(u,v)&=& \int_{-\infty}^{+\infty} \big{(}u^2+\frac23 bv^2
   \big{)}dx. \label{eq:cons2}
\end{eqnarray}
Later, Feng \cite{f:f} considered the
initial value problem for the following system:
\begin{equation}\label{eq:feng}
  \left \{
   \begin{array}{l}
     u_t-a(u_{xxx}+6uu_x)=2bvv_x, \\
     v_t+v_{xxx}+cuv_x+dvv_x=0, \\
     u(0)=u_0, \;\; v(0)=v_0,
   \end{array}
  \right.
\end{equation}
which reduces to the Hirota-Satsuma system when $c=3$ and $d=0$, always
assuming that $a \not =0$.  LWP of the
IVP associated to system (\ref{eq:feng}) was obtained,
for initial data $(u_0,v_0)\in H^s(\mathbb{R}) \times H^s(\mathbb{R})$
for $s \ge1$, with $a+1 \not = 0$ and $bc>0$.  Moreover, global
well-posedness (GWP) for system (\ref{eq:feng}) was also proved
(see \cite{f:f})  in $H^s(\mathbb{R})\times H^s(\mathbb{R})$
for $s \ge1$, if $-1<a<0$ and $bc>0$.

The second problem we will consider here is related to the local
well-posedness of the IVP for the Gear-Grimshaw system given by
\begin{equation}\label{eq:gg0}
  \left \{
   \begin{array}{l}
     u_t+u_{xxx}+a_{3}v_{xxx}+uu_x+a_1vv_x+a_2(uv)_{x}=0, \\
     b_1v_t+v_{xxx}+b_2a_3u_{xxx}+vv_x+b_2a_2uu_x+b_2a_1(uv)_x+rv_x=0,\\
     u(0)=u_0, \;\; v(0)=v_0,
   \end{array}
  \right.
\end{equation}
where $a_1, a_2, a_3 \in \mathbb{R}$, $r \in \mathbb{R}$, and
$b_1, b_2 >0$; $u=u(x,t)$, $v=v(x,t)$ are real-valued functions of
the two real variables $x$ and $t$.  System (\ref{eq:gg0}) was
derived in \cite{gg:gg} (see also \cite{bpst:bpst} for a very good explanation 
about the physical context in which this system arises) as a model to 
describe the strong interaction of two-dimensional, weakly nonlinear, 
long, internal gravity waves propagating on neighboring pycnoclines in 
a stratified fluid, where the two waves correspond 
to different modes.  Bona \textit{et al.} \cite{bpst:bpst} proved GWP 
of the IVP associated to (\ref{eq:gg0}) with initial
data belonging to $H^s(\mathbb{R}) \times H^s(\mathbb{R})$ for 
$s \ge 1$, assuming $r=0$ and $|a_3|<1/\sqrt{b_2}$.  Later, Ash 
\textit{et al.} \cite{acw:acw} considered GWP of (\ref{eq:gg0}) in 
$L^2(\mathbb{R}) \times L^2(\mathbb{R})$ supposing $r=0$ and 
$|a_3| \not = 1/\sqrt{b_2}$ (see Section \ref{subsection:comment}-{\bf{(2)}}). 
Further, Saut and Tzvetkov \cite{STz:STz} considered GWP of system 
(\ref{eq:gg}) for initial data $(u_0,v_0) \in L^2(\mathbb{R}) \times L^2(\mathbb{R})$, 
assuming that $r \not = 0$ and that the matrix $(a_{ij})_{i,j \in\{1,2\}}$ has real 
distinct eigenvalues (see Section \ref{subsection:comment}-{\bf{(1)}}).  Recently, 
Linares and Panthee \cite{lp:lp}, by using the bilinear estimate of Kenig, Ponce, 
and Vega \cite{kpv2:kpv2},  showed LWP for system (\ref{eq:gg2}) with 
initial data $(u_0,v_0) \in H^s(\mathbb{R}) \times H^s(\mathbb{R})$ for $s>-3/4$  
(see Section \ref{subsection:comment}-{\bf{(2)}}, and Remark 
\ref{correction}-{\bf{ii.)}}).  Solutions of  (\ref{eq:gg0}) satisfy the 
following conservation laws:
\begin{equation*}
 \Phi_1(u) = \int_{-\infty}^{+\infty} udx, \,\,\,
 \Phi_2(v) = \int_{-\infty}^{+\infty} v dx, \,\,\,
 \Phi_3(u,v)  = \int_{-\infty}^{+\infty} (b_2u^2+b_1v^2)dx,
\end{equation*}
\begin{equation*}
 \Phi_4(u,v)= \int_{-\infty}^{+\infty}
 \Big{(}b_2u_x^2+v_x^2+2b_2a_3u_xv_x-b_2\frac{u^3}{3}
 -b_2a_2 u^2v-b_2a_1uv^2-\frac{v^3}{3}-rv^2\Big{)}dx.
\end{equation*}

We say that the IVP
\begin{equation*}
   \left \{
   \begin{array}{l}
     \partial_t \vec{u} (t)=F(t, \vec{u}(t)), \\
     \vec{u}(0)=\vec{u}_0
   \end{array}
  \right.
\end{equation*}
is \textit{locally well-posed} in $X$ (Banach space) if there
exist $T=T(\|\vec{u}_0\|_X)>0$ and a unique solution $\vec{u}(t)$
of the corresponding IVP such that \\
{\bf{i.)}} $\vec{u} \in C([-T,T];X) \cap Y_T = X_T$; \\
{\bf{ii.)}} the mapping data-solution $\vec{u}_0 \mapsto
\vec{u}(t)$, from $\{ \vec{v}_0 \in X ; \|\vec{v}_0\|_X \le M \}$
into $X_T$ is uniformly continuous for all $M>0$; i.e.
\begin{eqnarray*}
  &&\forall M>0, \forall \epsilon >0, \exists \delta=\delta(\epsilon, M)>0,
    \|\vec{u}_0-\vec{v}_0\|_X < \delta \;\;\text{then} \;\;
    \|\vec{u}-\vec{v}\|_{X_T} < \epsilon, \\
  &&\text{where} \;\; \|\vec{u}_0\|_X \le M \;\;\text{and}\;\;
    \|\vec{v}_0\|_X \le M.
\end{eqnarray*}
We say that the IVP  is \textit{globally well-posed}
in $X$ if the same properties hold for all
time $T>0$. If some hypothesis in the definition of local
well-posedness fails, we say that the IVP is \textit{ill-posed}.

This paper is organized as follows.  In Section 2 we use Banach's
fixed-point theorem in a suitable function space and the theory
obtained by Kenig, Ponce, and Vega, to prove LWP 
to system (\ref{eq:hs}), for any $a, b \in \mathbb{R}$, with initial data
in $H^s(\mathbb{R}) \times H^s(\mathbb{R})$ for $3/4<s\le1$.  We also
show that system (\ref{eq:hs}) with $a \not = 0$ is ill-posed in
$H^s(\mathbb{R}) \times H^{s'}(\mathbb{R})$ for $s\in [-1,-\frac34)$, and
$s' \in \mathbb{R}$.  We begin Section 3 with a few comments to
scale changes carried out previously concerning the Gear-Grimshaw system.   
Thus, we introduce some Bourgain-type spaces $X_{s,b}^a$
for $a \not =0$, and $s,b \in \mathbb{R}$.   
Moreover, we prove some new mixed-bilinear estimates involving the two 
Bourgain-type spaces $X_{s,b}^{1}$ and $X_{s,b}^{-1}$ corresponding to 
$\partial_t-\partial_x^3$ and $\partial_t+\partial_x^3$ respectively, 
to obtain LWP for the Gear-Grimshaw system (\ref{eq:gg}) 
with $r=0$, $a_{12}=a_{21}=0$, $a_{11} = - a_{22} \not =0$, and initial 
data in $H^s(\mathbb{R}) \times H^s(\mathbb{R})$ for $s>-3/4$ 
(see Theorem \ref{teox2} below).   We remark that these mixed-bilinear 
estimates (see Proposition \ref{prop5}) presented here are not an 
immediate consequence of the estimates proved by Kenig, Ponce, and Vega 
in \cite{kpv2:kpv2} (see Remark \ref{Remark:non-equiv} and 
Remark \ref{Remark:general-indices, nonequivalence}-{\bf{ii.)}}).  Finally, 
we notice that system (\ref{eq:hs}) is treated separately from system 
(\ref{eq:gg0}) because the nonlinearity in (\ref{eq:hs}) has the non-divergence 
form, while the one in (\ref{eq:gg0}) has the divergence form; a possible 
difficulty with regard to the LWP of (\ref{eq:hs}) in lower Sobolev indices
could be related to the obtention of a suitable bilinear 
estimate for the nonlinear term in the second equation of (\ref{eq:hs}).

%%%%%%%%%%%%%%%%%%%%%%%%%%%%%%%%%%%%%%%%%%%%%%%%%%%%%%%%%%%%
\smallskip
\noindent
{\bf{Notation:}}
\begin{itemize}
\advance\leftskip by -2.2em

\advance\itemindent by 0em
\item $\hat f = \F f$ : the Fourier transform of $f$
      ($\F ^{-1}$ : the inverse of the Fourier transform), where \;
      $\hat f (\xi) = \frac{1}{\sqrt{2\pi}} \int e^{-i\xi
      x}f(x)dx \,\,$ for $f \in L^1(\mathbb{R})$.
\item $\| \cdot \|_s$, $(\cdot,\cdot)_s$: the norm and
      the inner product respectively in $H^s(\mathbb{R})$ (Sobolev
      space of order $s$ of $L^2$ type), $s \in \mathbb{R}$.
      $\|f\|_s^2 \equiv \int (1+|\xi|^2)^s |\hat f(\xi)|^2d\xi$. \\
      $\|\cdot\|=\|\cdot\|_0$: the $L^2(\mathbb{R})$ norm.
      $(\cdot,\cdot)$ denotes the inner product on
      $L^2(\mathbb{R})$.
\item $B(X,Y)$: set of bounded linear operators on $X$ to
      $Y$. If $X=Y$ we write $B(X)$.  $\|\cdot\|_{B(X,Y)}$:
      the operator norm in $B(X,Y)$.
\item $L^p=\{f;f \text{ is measurable on } \mathbb{R}, \|f\|_{L^p}
      <\infty\}$,
      where $\|f\|_{L^p} =\big{(}\int |f(x)|^pdx\big{)}^{1/p}$
      if $1 \le p<+\infty$, and
      $\|f\|_{L^{\infty}}= \text{ess} \sup_{x \in \mathbb{R}} |f(x)|$,
      $f$ is an equivalence class.
\item $C(I;X)$ : set of continuous functions on the interval
      $I$  into the Banach space $X$.
\item $\|f\|_{L^q_T L^p_x} \equiv \big{(}\int_{-T}^T\|f(\cdot,t)\|_{L^p}^q
      dt \big{)}^{1/q}$, \;\;\;\; $\|f\|_{L^q_t L^p_x} \equiv
      \|f\|_{L^q_T L^p_x}$  if $T=+\infty$.
\item $\|f\|_{L^p_xL^q_T} \equiv \big{\|}\big{(} \int_{-T}^T|f(\cdot,t)|^q
      dt \big{)}^{1/q} \big{\|}_{L^p}$, \;\;\;\;  $\|f\|_{L^p_x L^q_t} \equiv
      \|f\|_{L^p_x L^q_T}$  if $T=+\infty$.
\item $\langle \xi \rangle \equiv 1+|\xi|$,  for $\xi \in \mathbb{R}$.
\item Let $A, B$ be two $n \times n$ matrices.
      $A \sim B \;$ iff $\; \exists \, T \in GL(n)$, $T^{-1}A T=B$.
\end{itemize}

%%%%%%%%%%%%%%%%%%%%%%%%%%%%%%%%%%%%%%%%%%%%%%%%%%%%%%%%%%%%%%%%%%%%%%
% Local Well-Posedness and Ill-Posedness to the  Hirota-Satsuma System
%%%%%%%%%%%%%%%%%%%%%%%%%%%%%%%%%%%%%%%%%%%%%%%%%%%%%%%%%%%%%%%%%%%%%%
\section{On the Hirota-Satsuma System}
%%%%%%%%%%%%%%%%%%%%%%%%%%%%%%%%%%%%%%%%%%%%%%%%%%%%%%%%%%%%
\subsection{Local Well-Posedness}
Let us denote by
\begin{equation}\label{eq:group}
  U_a(t)= e^{at \partial_x^3},   \;\;\;
  \widehat{U_a(t)\phi}(\xi)=e^{-iat\xi^3}\hat \phi(\xi)
  \;\; \text{for} \;\; \phi\in H^s(\mathbb{R}),
\end{equation}
the group associated with the linear part of the first equation
of system (\ref{eq:hs}).  We note that $U(t)\equiv U_{-1}(t)$ is the
group associated with the linear part of the KdV equation.
Next theorem proves LWP to system (\ref{eq:hs}) in suitable Sobolev spaces.

%%%%%%%%%%%%%%%%%%%%%%%%%%%%%%%%%%
\begin{Theorem}\label{Theorem:lwp}
Let $a \not = 0$ and $3/4<s\le1$.  Then for any
$u_0, v_0\in H^s(\mathbb{R})$, there exists $T=T(\|u_0\|_s,\|v_0\|_s)>0$
(with $T(\rho, \eta) \rightarrow \infty$
as $\rho \rightarrow 0, \eta \rightarrow 0$) and a unique solution
$(u,v)$ of problem (\ref{eq:hs}) such that
\begin{eqnarray}
  && u, v \in C([-T,T];H^s(\mathbb{R})), \label{eq:lwp1}\\
  && u_x, v_x \in L^4_TL^{\infty}_x,  \label{eq:lwp2}\\
  && D^s_x u_x, D^s_x v_x \in L^{\infty}_x L^2_T,  \label{eq:lwp3}\\
  && u, v \in L^2_x L^{\infty}_T,  \label{eq:lwp4}\\
  && u_x, v_x \in L^{\infty}_x L^2_T.  \label{eq:lwp5}
\end{eqnarray}
For any $T' \in (0,T)$ there exist neighborhoods $V$ of
$u_0$ in $H^s(\mathbb{R})$ and $V'$ of $v_0$ in
$H^s(\mathbb{R})$ such that the map $(\tilde u_0,\tilde v_0)
\mapsto (\tilde u,\tilde v)$ from $V \times V'$ into the class defined by
(\ref{eq:lwp1})-(\ref{eq:lwp5}) with $T'$ instead of $T$ is Lipschitz. \\
If $u_0, v_0\in H^r(\mathbb{R})$ with $r>s$, then the above results
hold with $r$ instead of $s$ in the same time interval. \\
Moreover, from the conservation laws (\ref{eq:cons1}) and
(\ref{eq:cons2}) we can choose $T=+\infty$
at least for $s=1$, for $a+1>0$ and $b>0$.
\end{Theorem}
\begin{proof}
Let $\frac34 <s\le1$. Given $r \in \mathbb{R}$ and $T>0$, let us define
\begin{eqnarray}\label{eq:norm}
 \Lambda_r^T(u) &\equiv& \max_{[-T,T]} \|u(t)\|_r
 +\|u_x\|_{L_T^4L_x^{\infty}} + \|D^r_xu_x\|_{L_x^{\infty}L_T^2}
 \nonumber \\
 &&+(1+T)^{-1/2}\|u\|_{L_x^2L_T^{\infty}}
   +\|u_x\|_{L_x^{\infty}L_T^2}.
\end{eqnarray}
Denote $\||(u,v)|\| \equiv \Lambda_s^T(u)+\Lambda_s^T(v)$.
We consider the space
\begin{equation*}
  X^T = \Big{\{} (u,v) \in C([-T,T];H^s(\mathbb{R})) \times
  C([-T,T];H^s(\mathbb{R})) ;
   \||(u,v)|\| < \infty \Big{\}}
\end{equation*}
and
$  X_M^T = \{ (u,v) \in X^T ; \||(u,v)|\| \le M \}.$
Let us write the integral equations associated to problem
(\ref{eq:hs})
\begin{equation*}
  \left \{
   \begin{array}{l}
   \Phi_1(u,v)(t)=U_a(t) u_0 + \int_0^tU_a(t-t')
   (6auu_x+2bvv_x)(t')dt', \\
   \Phi_2(u,v)(t)=U(t)v_0-3\int_0^tU(t-t') (uv_x)(t')dt'.
   \end{array}
  \right.
\end{equation*}
We will prove that $\Phi:X_M^T \mapsto X_M^T$, where
$\Phi(u,v)\equiv(\Phi_1(u,v),\Phi_2(u,v))$, is a contraction map for
suitably chosen $M$ and $T$.   We have the following inequalities:
\begin{eqnarray}\label{eq:norms}
  \|U_a(t)u_0\|_r &\le& c \|u_0\|_r
      \;\;\;\; \text{for} \;\; r\in\mathbb{R},\label{eq:norms1} \\
  \|D^r_x \partial_x U_a(t)u_0\|_{L^{\infty}_x L^2_T}
      &\le& \frac{c}{|a|^{1/2}}\|D^r_x u_0\| \;\;\;\; \text{for} \;\; r\in\mathbb{R},
      \label{eq:norms2} \\
  \|\partial_x U_a(t)u_0\|_{L^4_TL^{\infty}_x}
      &\le& \frac{c}{|a|^{1/4}}\|u_0\|_r \;\;\;\; \text{for} \;\; r\ge3/4,
      \label{eq:norms3} \\
  \|U_a(t)u_0\|_{L^2_x L^{\infty}_T}
      &\le&c_{(a,r)}(1+T)^{1/2} \|u_0\|_r  \;\;\;\; \text{for} \;\; r>3/4.
      \label{eq:norms4}
\end{eqnarray}
Expression (\ref{eq:norms1}) is a group property.  Inequality (\ref{eq:norms2})
is a consequence of Theorem 4.1 in \cite{kpv0:kpv0}.  Expression (\ref{eq:norms3})
follows from  Theorem 2.1 in \cite{kpv0:kpv0}.  Estimate (\ref{eq:norms4})
is obtained by using Proposition 2.4 in \cite{la:la}.  It follows from
(\ref{eq:norms1})-(\ref{eq:norms4}) that
$\Lambda_s^T(U_a(t)u_0) \le c \|u_0\|_s$. Let $(u,v) \in X_M^T$. Then
\begin{eqnarray}\label{eq:lwp6}
  \Lambda_s^T(\Phi_1(u,v))
    &\le& c\,\|u_0\|_s +c\!\!\int_0^T \|(uu_x)(\tau)\|d\tau
          +c\!\!\int_0^T\|D^s_x(uu_x) (\tau)\|d\tau \nonumber \\
    &   & +c\!\!\int_0^T \|(vv_x)(\tau)\|d\tau
          +c\!\!\int_0^T\|D^s_x(vv_x) (\tau)\|d\tau.
\end{eqnarray}
Choose  $M\equiv 4c(\|u_0\|_s+\|v_0\|_s)$.  It follows that
\begin{equation}\label{eq:lwp7}
  \|uu_x\|_{L^2_TL^2_x} \le \|u_x\|_{L_x^{\infty}L_T^2}
                            \|u\|_{L_x^2L_T^{\infty}}
                        \le M^2(1+T)^{1/2}.
\end{equation}
Now, by using Theorem A.12 in \cite{kpv1:kpv1} and H\"older's
inequality, it follows that
\begin{eqnarray}\label{eq:lwp8}
  \|D^s_x(u u_x)\|_{L^2_TL^2_x}
     &\le& c \,\big{\|}\|u_x\|_{L_x^{\infty}}
             \|D_x^s u\|_{L^2_x}\big{\|}_{L_T^2}
             +\|u\|_{L_x^2L_T^{\infty}}
             \|D^s_x u_x\|_{L^{\infty}_x L^2_T} \nonumber \\
     &\le& c\, T^{1/4}\|D_x^s u\|_{L^{\infty}_TL^2_x}
           \|u_x\|_{L^4_TL^{\infty}_x}
           +M^2(1+T)^{1/2} \nonumber \\
     &\le& c\, M^2(T^{1/4}+(1+T)^{1/2}).
\end{eqnarray}
By replacing (\ref{eq:lwp7}) and (\ref{eq:lwp8}) (and
similar estimates for $v$) into (\ref{eq:lwp6}) we obtain
\begin{equation}\label{eq:lwp9}
  \Lambda_s^T(\Phi_1(u,v)) \le \frac{M}{4}
  +c\,M^2T^{1/2}(T^{1/4}+(1+T)^{1/2}).
\end{equation}
By choosing $T>0$ small enough such that
$T^{1/2}(T^{1/4}+(1+T)^{1/2}) \le \frac{1}{4cM}$, it follows that
$\Lambda_s^T(\Phi_1(u,v))\le \frac{M}{2}$.  Similarly we have that
$\Lambda_s^T(\Phi_2(u,v))\le \frac{M}{2}$.  Then,
for $M>0$ and $T>0$ chosen as above, $\Phi$ is a well-defined map
from $X_M^T$ to itself.  Analogously, we prove that $\Phi$ is a contraction map.
The rest of the proof is similar to the proof of Theorem 2.1 in
\cite{kpv1:kpv1}.
\end{proof}

%%%%%%%%%%%%%%%%%%%%%%%%%%%%%%%%%%
\begin{Theorem}\label{Theorem:lwpa=0}
Let $a = 0$ and $3/4<s\le1$.  Then for any
$u_0, v_0 \in H^s(\mathbb{R})$, there exists
$T=T(\|u_0\|_s,\|v_0\|_s)>0$ (with $T(\rho, \eta) \rightarrow \infty$
as $\rho \rightarrow 0, \eta \rightarrow 0$) and a unique solution
$(u,v)$ of problem (\ref{eq:hs}) such that
\begin{eqnarray}
  &&u,v \in C([-T,T];H^s(\mathbb{R})),   \label{eq:lwp10}\\
  &&v_x \in L^4_TL^{\infty}_x,  \label{eq:lwp11}\\
  &&D^s_x v_x \in L^{\infty}_x L^2_T,  \label{eq:lwp12}\\
  &&u,v \in L^2_x L^{\infty}_T,  \label{eq:lwp13}\\
  &&v_x \in L^{\infty}_x L^2_T. \label{eq:lwp14}
\end{eqnarray}
For any $T' \in (0,T)$ there exist neighborhoods $V$ of
$u_0$ in $H^s(\mathbb{R})$ and $V'$ of $v_0$ in
$H^s(\mathbb{R})$ such that the map $(\tilde u_0,\tilde v_0)
\mapsto (\tilde u,\tilde v)$ from $V \times V'$ into the class defined by
(\ref{eq:lwp10})-(\ref{eq:lwp14}) with $T'$ instead of $T$ is Lipschitz. \\
If $u_0, v_0\in H^r(\mathbb{R})$ with $r>s$, then the above results
hold with $r$ instead of $s$ in the same time interval. \\
If $s=1$ and $b>0$, then we can choose $T=+\infty$.
\end{Theorem}
\begin{proof}
Let $3/4<s\le1$.  Let $\Lambda_s^T(\cdot)$ be the norm defined by
(\ref{eq:norm}).  Denote by
\begin{equation*}
  \tilde \Lambda_s^T (u) \equiv \max_{[-T,T]} \|u(t)\|_s
  +\|u\|_{L^2_x L^{\infty}_T},
\end{equation*}
and  $\||(u,v)|\| \equiv \tilde \Lambda_s^T(u)+\Lambda_s^T(v)$.
Let $X^T$ and $X_M^T$ be defined  as in the proof of
Theorem \ref{Theorem:lwp}.  Let us now consider
$\Phi(u,v) \equiv (\Phi_1(u,v),\Phi_2(u,v))$, where
\begin{equation*}
  \left \{
   \begin{array}{l}
   \Phi_1(u,v)(t)=u_0 + 2b\int_0^t(vv_x)(t')dt', \\
   \Phi_2(u,v)(t)=U(t)v_0-3\int_0^tU(t-t') (uv_x)(t')dt'.
   \end{array}
  \right.
\end{equation*}
Let $(u,v)\in X_M^T$.  Then
\begin{equation*}
  \Lambda_s^T(\Phi_2(u,v)) \le c\|v_0\|_s
  +c \int_0^T \|(uv_x)(\tau)\|d\tau
  +c \int_0^T \|D^s_x (uv_x)(\tau)\|d\tau.
\end{equation*} 
We see that 
\begin{equation*}
  \|uv_x\|_{L^2_TL^2_x} \le \|v_x\|_{L^{\infty}_x L^2_T} 
  \|u\|_{L^2_xL^{\infty}_T} \le M^2.
\end{equation*}
Now, using Theorem A.12 in \cite{kpv1:kpv1} and H\"older's
inequality, we get 
\begin{eqnarray*}
 \|D^s_x(uv_x) \|_{L^2_T L^2_x} &\le& 
  c \,\big{\|}\|v_x\|_{L_x^{\infty}} \|D_x^s u\|_{L^2_x}\big{\|}_{L_T^2}
  +\|u\|_{L_x^2L_T^{\infty}} \|D^s_x v_x\|_{L^{\infty}_x L^2_T} \\
  &\le& c\, M^2(1+T^{1/4}).
\end{eqnarray*}
By choosing $M \equiv 6c(\|u_0\|_s+\|v_0\|_s)$ and $T>0$ such that
$T^{1/2}(T^{1/4}+(1+T)^{1/2})\le \frac{1}{6cM}$, we obtain
$\Lambda_s^T(\Phi_2(u,v)) \le \frac{M}{3}$ and
$\max_{[-T,T]}\|\Phi_1(u,v)(t)\|_s \le \frac{M}{3}$.   Moreover,
\begin{eqnarray*}
  \|\Phi_1(u,v)\|_{L^2_xL^{\infty}_T}
    &\le&  \|u_0\|+2b \| \int_0^T |(vv_x)(\tau)|d\tau \|_{L^2_x}
    \le  \|u_0\|+c T^{1/2} \|vv_x\|_{L^2_TL^2_x} \\
    &\le& \frac{M}{6} +cT^{1/2}M^2 (1+T)^{1/2} \le \frac{M}{3}.
\end{eqnarray*}
Then $\||\Phi(u,v)|\| \le M$.  The rest of the proof is as
for Theorem \ref{Theorem:lwp}.
\end{proof}

%%%%%%%%%%%%%%%%%%%%%%%%%%%%%%%%%%%%%%%%%%%%%%%%%%%%%%%%%%%%

\begin{Remark}
In \cite{Sk:Sk}, Sakovich considered the following system:
\begin{equation}\label{eq:sk}
  \left \{
   \begin{array}{l}
     u_{xxx}+auu_x+bvu_{x}+cuv_x+dvv_x+mu_t+nv_t=0, \\
     v_{xxx}+euu_x+fvu_{x}+guv_x+hvv_x+pu_t+qv_t=0,\\
     u(0)=u_0, \;\; v(0)=v_0,
   \end{array}
  \right.
\end{equation}
where $mq \neq np$. This system can be written as
\begin{align}\label{eq:sk1}
 \left(\!\!
 \begin{array}{c}
  u_{xxx} \\
  v_{xxx} \end{array}\!\!\!
  \right)
  +A_0\left( \!\!\begin{array}{c}
  uu_x \\
  vu_x \end{array}\!\!\! \right)+A_1\left(\!\! \begin{array}{c}
  uv_x \\
  vv_x \end{array}\!\!\! \right)+A_2\left(\!\! \begin{array}{c}
  u_t \\
  v_t \end{array} \!\!\!\right)=0,
\end{align}
where
\begin{align*}
A_0=\left( \begin{array}{cc}
  a & b\\
  e & f \end{array} \right), \;A_1=\left( \begin{array}{cc}
  c & d \\
  g & h \end{array} \right), \;A_2=\left( \begin{array}{cc}
  m & n \\
  p & q \end{array} \right).
\end{align*}
Since $A_2$ is nonsingular, multiplying (\ref{eq:sk1}) by
$A_2^{-1}$, we get
\begin{align*}
\left(\!\! \begin{array}{cc}
  u_t \\
  v_t \end{array}\!\!\! \right)+A_2^{-1}\left(\!\!\begin{array}{cc}
  u_{xxx} \\
  v_{xxx} \end{array}\!\!\! \right)+
  A_2^{-1}A_0\left(\!\! \begin{array}{cc}
  uu_x \\
  vu_x \end{array}\!\!\! \right)+A_2^{-1}A_1\left(\!\! \begin{array}{cc}
  uv_x \\
  vv_x \end{array}\!\!\! \right)=0.
\end{align*}
If $P\in GL(2)$ is such that
$P^{-1}A_2^{-1}P=\text{diag}(a_0, a_1)$, where $a_0$
and $a_1$ are the eigenvalues of $A_2^{-1}$, by making $U=(u,v)^t=PV$,
we obtain a new system of Hirota-Satsuma type.  Therefore,
similar results to  Theorems \ref{Theorem:lwp} and
\ref{Theorem:lwpa=0} are also valid for this new system.
\end{Remark}
%%%%%%%%%%%%%%%%%%%%%%%%%%%%%%%%%%%%%%%%%%%%%%%%%%%%%%%%%%%%
% Ill-Posedness: Hirota-Satsuma System
%%%%%%%%%%%%%%%%%%%%%%%%%%%%%%%%%%%%%%%%%%%%%%%%%%%%%%%%%%%%
\subsection{Ill-Posedness to the Hirota-Satsuma System}
%%%%%%%%%%%%%%%%%%%%%%%%%%%%%%%%%%%%%%%%%%%%%%%%%%%%%%%%%%%%
Let us remark that if $u(x,t)$ and $v(x,t)$ are solutions of
(\ref{eq:hs}), then $\tilde{u}(x,t)=\lambda^2u(\lambda x,
\lambda^3 t)$ and $\tilde{v}(x,t)=\lambda^2v(\lambda x, \lambda^3
t)$ are also solutions of (\ref{eq:hs}).   This scaling argument
suggests that the Cauchy problem for the Hirota-Satsuma system
is locally well-posed in $H^s(\mathbb{R}) \times H^{s'}(\mathbb{R})$
for $s, s' >-\frac32$. 
It is not difficult to see that the IVP associated to the KdV equation
\begin{equation*}
\left \{
 \begin{array}{l}
  w_t+w_{xxx}+6ww_x=0, \\
  w(x,0)=w_0(x)
  \end{array}
 \right.
\end{equation*}
is equivalent to the IVP
\begin{equation}\label{eq:kdv}
\left \{
 \begin{array}{l}
  u_t-a(u_{xxx}+6uu_x)=0, \\
  u(x,0)=u_0(x)=w_0(-x),
  \end{array}
 \right.
\end{equation}
through the transformation $u(x,t)=w(-x,at)$, for $a \not =0$.  Note that
if $u$ is a solution of (\ref{eq:kdv}), then $(u,0)$ is a solution of
problem (\ref{eq:hs}) with initial data $(u_0,0)$.  Then, it follows from
the ill-posedness result for the KdV equation (see \cite{cct:cct}) that 
the mapping data-solution associated to the IVP (\ref{eq:hs}) with 
$a \not = 0$ is not uniformly continuous in 
$H^s(\mathbb{R})\times H^{s'}(\mathbb{R})$, for  $s\in[-1, -\frac34)$, 
and $s' \in \mathbb{R}$.
%%%%%%%%%%%%%%%%%%%%%%%%%%%%%%%%%%%%%%%%%%%%%%%%%%%%%%%%%%%%
%On the Gear-Grimshaw System
%%%%%%%%%%%%%%%%%%%%%%%%%%%%%%%%%%%%%%%%%%%%%%%%%%%%%%%%%%%%
\section{On the Gear-Grimshaw System}
%%%%%%%%%%%%%%%%%%%%%%%%%%%%%%%%%%%%%%%%%%%%%%%%%%%%%%%%%%%%
%%%%%%%%%%%%%%%%%%%%%%%%%%%%%%%%%%%%%%%%%%%%%%%%%%%%%%%%%%%%%%%%%%
%%%%%%%%%%%%%%%%%%%%%%%%%%%%%%%%%%%%%%%%%%%%%%%%%%%%%%%%%%%%%%%%%%
\subsection{Initial Comments} 
\label{subsection:comment} 
{\bf{(1) }}
We consider the Gear-Grimshaw system given by
\begin{equation}\label{eq:gg}
  \left \{
   \begin{array}{l}
     u_t+a_{11}u_{xxx}+a_{12}v_{xxx}+b_1(uv)_{x}+b_2uu_x+b_3vv_x=0, \\
     v_t+a_{21}u_{xxx}+a_{22}v_{xxx}+rv_x+b_4(uv)_x+b_5uu_x+b_6vv_x=0,\\
     u(0)=u_0, \;\; v(0)=v_0.
   \end{array}
  \right.
\end{equation}
Suppose $r \not = 0$. Let $A, B$ and $C(U)$ be the matrices 
(see \cite{STz:STz}) defined by
\begin{equation*}
  A=\left( \begin{array}{cc}
  a_{11} & a_{12} \\
  a_{21} & a_{22} \end{array} \right), \;\;
  B=\left( \begin{array}{cc}
  0 & 0 \\
  0 & r \end{array} \right), \;\;
  C(U)=\left( \begin{array}{cc}
  b_2u+b_1v & b_1u+b_3v \\
  b_5u+b_4v & b_4u+b_6v \end{array} \right),
\end{equation*}
where $U=(u,v)^t$.  Let $T \in GL(2)$ such that
$T^{-1}AT=\text{diag}(\alpha_+, \alpha_-)$, where $\alpha_+$ and  
$\alpha_-$ are the eigenvalues of $A$, and  
$\alpha_+, \alpha_- \in \mathbb{R}$.  By making 
$U=TV$, we obtain 
\begin{equation}\label{eq:well}
  \left \{
   \begin{array}{l}
     \!V_t(t,x)+\text{diag}(\alpha_+,\alpha_-)V_{xxx}(t,x)
     +B_1V_x(t,x)+C_1(V)(t,x)V_x(t,x)=0,  \\
     \!V(0)=T^{-1}U_0,
   \end{array}
  \right.
\end{equation}
where $B_1=T^{-1}BT = (b_{ij})_{i,j \in\{1,2\}}$, 
$C_1(V)=T^{-1}C(TV)T$ and $U_0=U(0)$.  Let $V=(v_1,v_2)^t$.  If 
we make the scale change (supposing 
$\alpha_+ \not =0, \alpha_- \not = 0$) 
\begin{equation*}
  v_1(t,x)=w_1\Big{(}t, \frac{x}{\alpha_+^{1/3}}\Big{)}, \;\;\;
  v_2(t,x)=w_2\Big{(}t, \frac{x}{\alpha_-^{1/3}}\Big{)},
\end{equation*}
then $W=(w_1,w_2)^t$ satisfies the following system:
\begin{equation*}\label{eq:correc}
  \left \{
   \begin{array}{l}
     \partial_1w_1(t,\frac{x}{\alpha_+^{1/3}})
     + \partial_2^3 w_1(t,\frac{x}{\alpha_+^{1/3}})
     + \frac{b_{11}}{\alpha_+^{1/3}} \partial_2w_1(t,\frac{x}{\alpha_+^{1/3}})
     + \frac{b_{12}}{\alpha_-^{1/3}} \partial_2w_2(t,\frac{x}{\alpha_-^{1/3}})
     + ...=0,\\
     \partial_1w_2(t,\frac{x}{\alpha_-^{1/3}})
     + \partial_2^3 w_2(t,\frac{x}{\alpha_-^{1/3}})
     + \frac{b_{21}}{\alpha_+^{1/3}} \partial_2w_1(t,\frac{x}{\alpha_+^{1/3}})
     + \frac{b_{22}}{\alpha_-^{1/3}} \partial_2w_2(t,\frac{x}{\alpha_-^{1/3}})
     + ...=0,
   \end{array}
  \right.
\end{equation*}
where $\partial_i$, for $i=1,2$ denotes the partial derivative
with respect to the $i$-th variable.  
It should be noted that $\partial_2w_1$ is evaluated at the point 
$(t,\frac{x}{\alpha_+^{1/3}})$ and $\partial_2w_2$ is evaluated 
at the point $(t,\frac{x}{\alpha_-^{1/3}})$.  Take 
$b_1=...=b_6=0$ in (\ref{eq:gg}).  
If $\alpha_+ \not = \alpha_-$, it follows that 
we should take care in any of the following cases:
\begin{itemize}
\item $b_{12} \not =0$ and $\partial_2w_2(t,\frac{x}{\alpha_+^{1/3}})
\not =\partial_2w_2(t,\frac{x}{\alpha_-^{1/3}})$, \\
\item $b_{21}\not =0$ and $\partial_2w_1(t,\frac{x}{\alpha_+^{1/3}})
\not = \partial_2w_1(t,\frac{x}{\alpha_-^{1/3}})$.
\end{itemize}

\noindent
{\bf{(2) }}
We now consider the following system
($C(U) \not \equiv 0$ and $r=0$ in (\ref{eq:gg})):
\begin{equation}\label{eq:gg1}
  \left \{
   \begin{array}{l}
     u_t+u_{xxx}+a_{3}v_{xxx}+uu_x+a_1vv_x+a_2(uv)_{x}=0, \\
     b_1v_t+v_{xxx}+b_2a_3u_{xxx}+vv_x+b_2a_2uu_x+b_2a_1(uv)_x=0,\\
     u(0)=u_0, \;\; v(0)=v_0,
   \end{array}
  \right.
\end{equation}
where $a_1, a_2, a_3, b_1, b_2$  are real constants, with $b_1, b_2>0$,
$a_3 \not = 0$, and $a_3^2b_2 \not = 1$.  We define 
(see \cite{acw:acw} and \cite{lp:lp}): 
  $\lambda = \{ ( 1 -\frac{1}{b_1} )^2
  +\frac{4b_2a_3^2}{b_1} \}^{1/2}$ 
  and
  $\alpha_{\pm} = \frac12 (1+\frac{1}{b_1} \pm\lambda)$.
Consider
\begin{equation*}
  \left \{
   \begin{array}{l}
     \tilde u(t,x)=(\frac{1-\alpha_-}{\lambda})u (t,\alpha_+^{1/3}x)
     +\frac{a_3}{\lambda}v(t,\alpha_+^{1/3}x), \\
     \tilde v(t,x)=(\frac{\alpha_+-1}{\lambda}) u (t,\alpha_-^{1/3}x)
     -\frac{a_3}{\lambda} v(t, \alpha_-^{1/3}x),
   \end{array}
  \right.
\end{equation*}
or equivalently
\begin{equation}\label{eq:scale}
  \left \{
   \begin{array}{l}
     u(t,x)=\tilde u (t,\frac{x}{\alpha_+^{1/3}})
     +\tilde v(t,\frac{x}{\alpha_-^{1/3}}), \\
     v(t,x)=(\frac{\alpha_+-1}{a_3}) \tilde u (t,\frac{x}{\alpha_+^{1/3}})
     -(\frac{1-\alpha_-}{a_3}) \tilde v(t, \frac{x}{\alpha_-^{1/3}}).
   \end{array}
  \right.
\end{equation}
We note that this change of variable is equivalent to the one  
performed in item  {\bf{(1)}} for $W$.  
Take $b_1=b_2=1$, $a_1=a_2=0$ and $a_3=2$ in system
(\ref{eq:gg1}).  Then $\alpha_+=3$, $\alpha_-=-1$ and $\lambda=4$.
By using (\ref{eq:scale}), it follows that
\begin{equation*}
  \left \{
   \begin{array}{l}
     \partial_1\tilde u_ (t,\frac{x}{3^{1/3}})+\partial_1\tilde v_(t,-x)
     +\partial_2^3\tilde u (t,\frac{x}{3^{1/3}})+\partial_2^3\tilde v (t,-x)
     +\frac{1}{3^{1/3}} \tilde u(t,\frac{x}{3^{1/3}})
     \partial_2 \tilde u(t,\frac{x}{3^{1/3}}) \\
     -\tilde u(t,\frac{x}{3^{1/3}}) \partial_2 \tilde v (t,-x)
     +\frac{1}{3^{1/3}} \partial_2 \tilde u(t,\frac{x}{3^{1/3}})\tilde v (t,-x)
     -\tilde v(t,-x) \partial_2 \tilde v (t,-x)=0,\\
     \partial_1\tilde u_ (t,\frac{x}{3^{1/3}})-\partial_1\tilde v_(t,-x)
     +\partial_2^3\tilde u (t,\frac{x}{3^{1/3}})-\partial_2^3\tilde v (t,-x)
     +\frac{1}{3^{1/3}} \tilde u(t,\frac{x}{3^{1/3}})
     \partial_2 \tilde u(t,\frac{x}{3^{1/3}}) \\
     +\tilde u(t,\frac{x}{3^{1/3}}) \partial_2 \tilde v (t,-x)
     -\frac{1}{3^{1/3}} \partial_2 \tilde u(t,\frac{x}{3^{1/3}})\tilde v (t,-x)
     -\tilde v(t,-x) \partial_2 \tilde v (t,-x)=0.
   \end{array}
  \right.
\end{equation*}
Then
\begin{equation*}\label{eq:correcpanthee}
   \left \{
   \begin{array}{l}
     \partial_1\tilde u_ (t,\frac{x}{3^{1/3}})
     +\partial_2^3\tilde u (t,\frac{x}{3^{1/3}})
     +\frac{1}{3^{1/3}} \tilde u(t,\frac{x}{3^{1/3}})
     \partial_2 \tilde u(t,\frac{x}{3^{1/3}})
     -\tilde v(t,-x) \partial_2 \tilde v (t,-x)=0,\\
     \partial_1\tilde v_(t,-x)+\partial_2^3\tilde v (t,-x)
     -\tilde u(t,\frac{x}{3^{1/3}}) \partial_2 \tilde v (t,-x)
     +\frac{1}{3^{1/3}} \partial_2 \tilde u(t,\frac{x}{3^{1/3}})
     \tilde v (t,-x)=0,\\
     \tilde u(0,x)=\frac12 u_0(3^{1/3}x) +\frac12 v_0(3^{1/3}x), \\
     \tilde v(0,x)=\frac12 u_0(-x)- \frac12 v_0(-x).
   \end{array}
  \right.
\end{equation*}
Notice that 
 $-\tilde u(t,\frac{x}{3^{1/3}}) \partial_2 \tilde v (t,-x)
 +\frac{1}{3^{1/3}} \partial_2 \tilde u(t,\frac{x}{3^{1/3}})
 \tilde v (t,-x)=\partial_x(\tilde u(t,\frac{x}{3^{1/3}}) 
 \tilde v (t,-x)),$
where $\partial_x \neq \partial_2$.   It follows that, in general,  
system (\ref{eq:gg1}) cannot be written as
\begin{equation}\label{eq:gg2}
  \left \{
   \begin{array}{l}
     \tilde u_t+\tilde u_{xxx}+a\tilde u\tilde u_x+b\tilde v\tilde v_x
     +c(\tilde u \tilde v)_x=0, \\
     \tilde v_t+\tilde v_{xxx}+\tilde a \tilde u \tilde u_x
     +\tilde b \tilde v \tilde v_x+\tilde c (\tilde u \tilde v)_x=0,\\
     \tilde u(0,x)=(\frac{1-\alpha_-}{\lambda})u_0(\alpha_+^{1/3}x)
     +\frac{a_3}{\lambda}v_0(\alpha_+^{1/3}x), \\
     \tilde v(0,x)=(\frac{\alpha_+-1}{\lambda})u_0(\alpha_-^{1/3}x)
     -\frac{a_3}{\lambda}v_0(\alpha_-^{1/3}x),
    \end{array}
  \right.
\end{equation}
where $a, b, c$ and $\tilde a, \tilde b, \tilde c$ are
constants.  

\begin{Remark}\label{correction}
{\bf{i.)}} To prove LWP to a system like (\ref{eq:gg}) with $r=0$, we 
can work with an equivalent system like (\ref{eq:well}) (see Remark 
\ref{Remark:No-mixedterms}). In this case 
and if $\alpha_+, \alpha_- \in \mathbb{R} \setminus \{0\}$,  we can  
consider the two groups $U_{-\alpha_+}(t)=e^{-(\alpha_+)t\partial_x^3}$
and $U_{-\alpha_-}(t)=e^{-(\alpha_-)t\partial_x^3}$ associated to
the linear part of system (\ref{eq:well}) (see Theorem \ref{teox1} 
and Corollary \ref{Corollary:teox2} for the case when 
$|\alpha_+|=|\alpha_-|$). \\  
{\bf{ii.)}} The LWP result obtained in \cite{lp:lp} really corresponds 
to system (\ref{eq:gg2}).  To prove LWP for the more general case 
corresponding to system (\ref{eq:gg}) with $r=0$, we could try to obtain  
some suitable bilinear estimates (see Propositions \ref{prop4} and 
\ref{prop5}, and 
Remark \ref{Remark:general-indices, nonequivalence}-{\bf{i.)}} for 
the case when $|\alpha_+|=|\alpha_-| \not=0$).
\end{Remark}

%%%%%%%%%%%%%%%%%%%%%%%%%%%%%%%%%%%%%%%%%%%%%%%%%%%%%%%%%%%%%%%%%%
%Some Properties of  X_{s,b}^a-Spaces
%%%%%%%%%%%%%%%%%%%%%%%%%%%%%%%%%%%%%%%%%%%%%%%%%%%%%%%%%%%%%%%%%%
\subsection{Definiton of $X_{s,b}^a$-Spaces}
Let $a \not =0$.  For $s,b\in \mathbb{R}$, $X_{s,b}^a$ is used to
denote the completion of the Schwartz space
$\mathcal{S}(\mathbb{R}^2)$ with respect to the norm
\begin{equation}\label{eq:norm1}
  \|F\|_{X_{s,b}^a} \equiv \Big{(} \int_{-\infty}^{+\infty}
  \int_{-\infty}^{+\infty} \langle \tau+a\xi^3 \rangle^{2b}
  \langle \xi \rangle^{2s}  |\widehat F(\xi,\tau)|^2 d\xi d\tau
  \Big{)}^{1/2},
\end{equation}
where $\widehat F(\xi,\tau)=(2\pi)^{-1} \int_{\mathbb{R}^2}
e^{-i(x\xi+t\tau)}F(x,t)dxdt$.  It follows that $X_{s,b}^{-1}$
coincides with the usual Bourgain space $X_{s,b}$ for the KdV
equation (see \cite{bou:bou}).

%%%%%%%%%%%%%%%%%%%%%%%%%%%%%%%%%%%%%%%%%%%%%%%%%%%%%

\begin{Lemma}\label{L1}
Let $b>1/2$, $s\geq-3/2$, and $a_0, a_1 \in \mathbb{R} \setminus \{0\}$ 
such that $a_0\neq a_1$.  Then
\begin{equation*}
  X_{s,b}^{a_0} \neq X_{s,b}^{a_1}.
\end{equation*}
\end{Lemma}
\begin{proof}
First, we suppose that $a_0 \cdot a_1<0$.  We may assume that $a_0>0$. \\
{\bf Case: $s>1/2-b$.} Consider
$v \in X_{s,b}^{a_1} \cap L^2(\mathbb{R}^2)$ such that
\begin{equation*}
 |\hat{v}(\xi,\tau)|^{2}=\frac{1}{\langle \xi
 \rangle^{2s+2b}\langle \tau+a_1\xi^3 \rangle^{4b}}.
\end{equation*}
Therefore
\begin{equation*}
  \|v\|_{X_{s,b}^{a_0}}^2\geq c(b,a_0) \int_{(\mathbb{R}^+)^2}\frac{\xi^{6b}
  d\xi d\tau}{\langle \xi \rangle^{2b} \langle \tau+a_1\xi^3
  \rangle^{4b}}=\infty.
\end{equation*}
{\bf Case: $-3/2 \le s\le 0$.}
Consider $u \in X_{s,b}^{a_1}\cap L^2(\mathbb{R}^2)$ such that
\begin{equation*}
 |\hat{u}(\xi,\tau)|^{2}=\frac{1}{\langle \xi
 \rangle^{d} \langle \tau+a_1\xi^3 \rangle^{4b}}, \;\;\;
 d \in (1,6b-2).
\end{equation*}
Therefore
\begin{equation*}
  \|u\|_{X_{s,b}^{a_0}}^2\geq c(b,a_0)
  \int_{(\mathbb{R}^+)^2}\frac{\xi^{6b} \langle \xi \rangle^{2s-d}
  d\xi d\tau}{\langle \tau+a_1\xi^3
  \rangle^{4b}}=\infty.
\end{equation*}
The case $a_0 \cdot a_1>0$ follows from the case $a_0 \cdot a_1<0$  
and from Lemma \ref{L2}.
\end{proof}
%%%%%%%%%%%%%%%%%%%%%%%%%%%%%%%%%%%%%%%%%%%%%%%%%
\begin{Remark}\label{Remark:non-equiv}
Lemma \ref{L1} implies that the two norms 
$\| \cdot\|_{X^{1}_{s,b}}$ and $\| \cdot\|_{X^{-1}_{s,b}}$ are 
not equivalent for $s>-3/4$ and $b>1/2$.  Then, it follows that  
Proposition \ref{prop5} below is not an immediate consequence of Proposition 
\ref{prop4}.
\end{Remark}

%%%%%%%%%%%%%%%%%%%%%%%%%%%%%%%%%%%%%%%%%%%%%%%%%%%%%%%%%%%%%%%%%%
%Bilinear Estimates in  X_{s,b}^a-Spaces
%%%%%%%%%%%%%%%%%%%%%%%%%%%%%%%%%%%%%%%%%%%%%%%%%%%%%%%%%%%%%%%%%%
\subsection{Bilinear Estimates in $X_{s,b}^a$-Spaces}
%%%%%%%%%%%%%%%%%%%%%%%%%%%%%
\begin{Proposition}\label{prop4}
Given $s>-\frac{3}{4}$ and $a\neq 0$, there exist
$b' \in (-\frac{1}{2},0)$ and $\epsilon_s
>0$ such that for any $b \in(\frac{1}{2}, b'+1]$ with $b'+1-b\le \epsilon_s$ 
\begin{equation}\label{eq:bil0}
 \|(uv)_x\|_{X_{s,b'}^a}\quad \leq \quad c_{(a,s,b,b')}\,\|u\|_{X_{s,b}^a} 
 \,\|v\|_{X_{s,b}^a}.
\end{equation}
\end{Proposition}

\begin{proof}
The result follows from Corollary 2.7 in \cite{kpv2:kpv2}, and 
from the fact that if $g(x,t)\equiv f(x,\frac{t}{-a})$ then 
$\widehat{g}( \xi,\tau)=|a|\widehat{f}(\xi,-a\tau)$.
\end{proof}

The following lemma contains elementary calculus inequalities.
%%%%%%%%%%%%%%%%%%%%%%%%
\begin{Lemma}\label{Lemma:Calculus}
If $b>1/2$, then there exists $c_b>0$ such that
\begin{equation}\label{eq:Calculus1}
  \int \frac{dx}{(1+|a| |x^2-\eta^2|)^{2b}} \le \frac{c_b}{|a| |\eta|}.
\end{equation}
If $0 \le \alpha \le \beta$, and $\beta>1$, then there exists
$c_{(\alpha,\beta)}>0$ such that
\begin{equation}\label{eq:Calculus2}
  \int \frac{dx}{(1+|x-a'|)^{\alpha} (1+|x-a|)^{\beta}} \le 
  \frac{c_{(\alpha,\beta)}}{(1+|a-a'|)^{\alpha}}.
\end{equation}
\end{Lemma}
\begin{proof}
To prove (\ref{eq:Calculus1}) we consider the two integrals
corresponding to $|x-\eta|>|\eta|$ and $|x-\eta|\le |\eta|$.  To prove
(\ref{eq:Calculus2}) we may suppose $a'=0$, then we consider the 
integrals corresponding to $|x|>|a|/2$ and $|x| \le |a|/2$ (see (2.12) in 
\cite{bop:bop}).
\end{proof}

Next lemma will be useful for the proof of Lemma \ref{L4}.
%%%%%%%%%%%%%%%%%%%%%%%%

\begin{Lemma}\label{Lemma:Estimative-L4}
If $s\in [-\frac34,0]$, $b' \le \frac{s}{3} -\frac14$ and $b>\frac12$, then there
exists $c_{(s,b,b')}>0$ such that
\begin{equation}\label{eq:Estimative-L4}
 \phi_1(\xi,y) \equiv 
  \frac{|\xi|^{3-4s}}{\langle \xi^3(y+2) \rangle^{-2b'} \langle \xi 
  \rangle^{-2s}}
  \int \frac{|y+2|^{-2s}dx}{\langle \xi^3(y+3/4-x^2) \rangle^{2b}} \le
  c_{(s,b,b')}.
\end{equation}
\end{Lemma}
\begin{proof}
{\bf{i.)}} First, we suppose $|y+\frac34|>\frac14$.  Since $s \le0$ and 
$b\ge0$, it follows that
\begin{equation*}
  \phi_1(\xi,y) \le
  \frac{|\xi|^{3-4s} |y+2|^{-2s}}{\langle \xi^3(y+2) \rangle^{-2b'}} 
  \int \frac{dx}{(1+ |\xi|^3|x^2-(|y+3/4|^{1/2})^2|)^{2b}}.
\end{equation*}
Since $b>1/2$, it follows from (\ref{eq:Calculus1}) that 
\begin{equation*}
  \phi_1(\xi,y) \le
  c_b \frac{(|\xi|^3|y+2|)^{-4s/3} |y+2|^{-2s/3}}{\langle 
  \xi^3(y+2) \rangle^{-2b'}|y+3/4|^{1/2}}
  \le c_{(s,b,b')} \frac{|y+2|^{-2s/3}}{|y+3/4|^{1/2}},
\end{equation*}
where in the last inequality we have used 
$b' \le \frac{s}{3}-\frac14$ and $s\ge-\frac34$. \\
The case $|y+2|\le 5/2$ is immediate.  If 
$|y+2|>5/2$, then $|y+3/4|\ge |y+2|-5/4>|y+2|/2$; hence 
$\phi_1(\xi,y)\le c_{(s,b,b')}|y+2|^{-\frac{2s}{3}-\frac12} 
\le c_{(s,b,b')}$, for $s \ge-3/4$. \\  
{\bf{ii.)}} Second, we suppose $-\frac14 \le y+\frac34 \le 0$.  Since 
$s \le0$ and $b>\frac14$, it follows that 
\begin{eqnarray*}
  \phi_1(\xi,y) &\le& c_b 
  \frac{|\xi|^{\frac32-2s}|y+2|^{-2s}}{\langle 
  \xi^3(y+2) \rangle^{-2b'}}
  \int_0^{+\infty}\!\!\! \frac{|\xi|^{\frac32}dx}
  {(1+|\xi|^3(|y+\frac34|^{\frac12}+x)^2)^{2b}} \\
  &\le&  c_b 
  \frac{|\xi|^{\frac32-2s}|y+2|^{-2s}}{\langle 
  \xi^3(y+2) \rangle^{-2b'}}
  \int_0^{+\infty} \!\!\!\frac{dz}{(1+z^2)^{2b}} 
  \le c_b \frac{|\xi|^{\frac32-2s}|y+2|^{-2s}}{\langle 
  \xi^3(y+2) \rangle^{-2b'}}.
\end{eqnarray*}
Since $1\le y+2 \le \frac54$, $b'\le 0$, and $s\le0$, it follows 
from the last inequality that 
\begin{equation*}
  \phi_1(\xi,y) \le c_{(s,b)} \frac{|\xi|^{\frac32-2s}}
  {\langle \xi^3 \rangle^{-2b'}} \le c_{(s,b,b')},
\end{equation*}
where the last inequality is a consequence of the fact that 
$b'\le\frac{s}{3}-\frac14$. \\
{\bf{iii.)}} Finally, we consider the case 
$0< y+\frac34 \le \frac14$.  Since $\frac54 <y+2 \le\frac32$, 
$s\le0$ and $b'\le \frac{s}{3}-\frac14$, it follows that
\begin{eqnarray*}
  \phi_1(\xi,y) &\le& c_s \frac{|\xi|^{\frac32-2s}}
  {(1+|\xi|^3)^{-2b'}} \int_0^{+\infty} \frac{|\xi|^{\frac32}dx}
  {(1+|\xi|^3(y+\frac34)|1-\frac{x^2}{y+3/4}|)^{2b}} \\
  &\le& c_{(s,b,b')}
  \int_0^{+\infty} \frac{|\xi|^{\frac32}(y+\frac34)^{\frac12}dz}
  {(1+|\xi|^3(y+\frac34)|1-z^2|)^{2b}}.
\end{eqnarray*}
Now, we split the last integral into two parts, namely $|z|\le \sqrt 2$ 
and $|z|> \sqrt 2$.  Since $2b>\frac12$, it follows that 
\begin{equation*}
  \int_0^{\sqrt 2} \frac{|\xi|^{\frac32}(y+\frac34)^{\frac12}dz}
  {(1+|\xi|^3(y+\frac34)|1-z^2|)^{2b}} 
  \le \int_0^{\sqrt 2} \frac{dz}{|1-z^2|^{1/2}}
  \le c \int_0^{\sqrt 2} \frac{dz}{|1-z|^{1/2}} \le c.
\end{equation*}
On the other hand, since $z^2>2$ implies $z^2-1>z^2/2$, 
and by making the change of variable 
$x=|\xi|^{\frac32}(y+\frac34)^{\frac12}z$, it follows that
\begin{equation*}
  \int_{\sqrt 2}^{+\infty} \frac{|\xi|^{\frac32}(y+\frac34)^{\frac12}dz}
  {(1+|\xi|^3(y+\frac34)|1-z^2|)^{2b}} 
  \le c_b\int_0^{+\infty} \frac{dx}{(1+x^2)^{2b}} =c_b.
\end{equation*}
\end{proof}

%%%%%%%%%%%%%%%%%%%%%%%%

The next eight lemmas will be used for proving Proposition
\ref{prop5}. 
% Their statements are analogous to Lemmas 2.4-2.6 in \cite{kpv2:kpv2}.

%%%%%%%%%%%%%%%%%%%%%%%%
\begin{Lemma}\label{L3}
If $b' \le -\frac14$ and $b>\frac12$, then there exists $c_b>0$ such that
\begin{equation}\label{eq:L3-0}
 \frac{|\xi|}{\langle \tau+\xi^3 \rangle^{-b'}} \Big(  \int \!\!\! \int 
 \frac{d\xi_1 d\tau_1}{\langle \tau_1-\xi_1^3 \rangle^{2b} 
  \langle \tau-\tau_1-(\xi-\xi_1)^3 \rangle^{2b}}\Big)^{1/2} \le c_b.
\end{equation}
\end{Lemma}
\begin{proof}
Since $b>1/2$, it follows from (\ref{eq:Calculus2}) that
\begin{equation*}
  \int \frac{d \tau_1}{\langle \tau_1-\xi_1^3 \rangle^{2b} 
 \langle \tau-\tau_1-(\xi-\xi_1)^3 \rangle ^{2b}} \le
\frac{c_b}{\langle \tau-\xi^3+3\xi \xi_1(\xi-\xi_1) \rangle^{2b}}.
\end{equation*}
Then, it is sufficient to prove that
\begin{equation*}
  \frac{|\xi|^2}{(1+|\tau+\xi^3|)^{-2b'}} 
  \int \frac{d \xi_1}{(1+|\tau-\xi^3+3\xi \xi_1(\xi-\xi_1)|)^{2b}} \le c.
\end{equation*}
By making the change of variable $\tau=\xi^3(1+z)$, 
we need now to verify the following:
\begin{equation*}
  \frac{|\xi|^2}{(1+|\xi^3||z+2|)^{-2b'}} 
  \int \frac{d \xi_1}{(1+|\xi^3 z+3\xi \xi_1(\xi-\xi_1)|)^{2b}} \le c.
\end{equation*}
By performing the change of variable $\xi_1=\xi x$ inside the last 
integral, and $z=3y$,  and since  $x-x^2=\frac14 - (x-\frac12)^2$, 
it is not difficult to see that the expression
we need to prove now is the following
\begin{equation*}
 \phi(\xi,y)\equiv \frac{|\xi|^3}{(1+|\xi|^3|3y+2|)^{-2b'}}
 \int \frac{dx}{(1+|\xi|^3|y+1/4-x^2|)^{2b}} \le c.
\end{equation*}
{\bf{i.)}} First, we consider the case $|y+1/4|>1/12$.   Then
\begin{equation*}
  \phi(\xi,y)\le \frac{|\xi|^3}{(1+|\xi|^3|3y+2|)^{-2b'}}
 \int \frac{dx}{(1+|\xi|^3|x^2-(|y+1/4|^{1/2})^2|)^{2b}} \le c_b,
\end{equation*}
where in the last inequality we have used (\ref{eq:Calculus1}) and $b'
\le 0$.  \\
{\bf{ii.)}} Second, we assume  $-1/3 \le  y\le -1/4$.  Since $1+|\xi|^3
|3y+2|\ge 1+|\xi|^3$ and $b'\le -1/4$, it follows that 
$\frac{|\xi|^{3/2}}{(1+|\xi|^3|3y+2|)^{-2b'}} \le \frac{|\xi|^{3/2}}
{(1+|\xi|^3)^{-2b'}} \le 1$.  Then   
\begin{equation*}
  \phi(\xi,y)\le c_b \int_0^{+\infty} \frac{|\xi|^{3/2} dx}
  {(1+|\xi|^3 (|y+1/4|^{1/2}+x)^2)^{2b}} 
  \le c_b \int_0^{+\infty} \frac{dz}{(1+z^2)^{2b}} \le c_b.
\end{equation*}
{\bf{iii.)}} Finally, we suppose $-1/4 <y \le -1/6$.  Since $b' \le-1/4$,
and by making the change of variable $x=(y+1/4)^{1/2}z$, we get  
\begin{equation*}
  \phi(\xi,y)\le c \int_0^{+\infty} 
  \frac{|\xi|^{3/2}(y+1/4)^{1/2}dz}{(1+|\xi|^3(y+1/4) |1-z^2|)^{2b}}
  \le c,
\end{equation*}
where in the last inequality we have used the following estimates.  
Since $b\ge 1/4$, it follows that  $(|\xi|^3(y+1/4)|1-z^2|)^{1/2} \le
(1+|\xi|^3(y+1/4)|1-z^2|)^{2b}$.  Then 
\begin{equation*}
 \int_0^{\sqrt 2}  \frac{|\xi|^{3/2}(y+1/4)^{1/2}dz}{(1+|\xi|^3(y+1/4)
 |1-z^2|)^{2b}} \le \int_0^{\sqrt 2} \frac{dz}{|1-z|^{1/2}|1+z|^{1/2}}
 \le c.
\end{equation*}
Moreover, since $z^2-1>z^2/2$ for $z>\sqrt 2$, and $b>1/2$, it follows that
\begin{equation*}
 \int_{\sqrt 2}^{+\infty} \frac{|\xi|^{3/2}(y+1/4)^{1/2}dz}{(1+|\xi|^3(y+1/4)
 |1-z^2|)^{2b}} \le \int_{\sqrt 2}^{+\infty}
 \frac{2|\xi|^{3/2}(y+1/4)^{1/2}dz}
 {1+|\xi|^3(y+1/4)z^2}\le \int_0^{+\infty} \!\!\! \frac{2dx}{1+x^2}.
\end{equation*}

\end{proof}

%%%%%%%%%%%%%%%%%%%%%%%%
\begin{Lemma}\label{L4}
If $s \in[-\frac34,-\frac14]$, $b' \in [-\frac12,\frac{s}{3}-\frac14]$  
and $b>\frac12$, then there exists $c_{(s,b,b')}>0$ such that
\begin{equation}\label{eq:L4-0}
 \frac{|\xi|}{\langle  \tau + \xi^3 \rangle^{-b'} \langle \xi
 \rangle^{-s}}  \Big( \int\!\!\!\int_A 
 \frac{ |\xi_1(\xi-\xi_1)|^{-2s}d\tau_1 d\xi_1}{\langle \tau_1-\xi_1^3
 \rangle^{2b}  \langle \tau-\tau_1-(\xi-\xi_1)^3 \rangle^{2b}}\Big)^{1/2}
 \le c_{(s,b,b')},
\end{equation}
where $A=A(\xi,\tau)$ is defined as
\begin{equation*}
  A=\{ (\xi_1,\tau_1) \in \mathbb{R}^2 ; |\xi_1| \ge 1, |\xi-\xi_1|
  \ge 1, |\tau-\tau_1-(\xi-\xi_1)^3| \le |\tau_1-\xi_1^3| \le
  |\tau+\xi^3| \}.
\end{equation*}
\end{Lemma}
\begin{proof}
We denote by $\chi_D$ the characteristic function of the
set $D$.  We remark that  $A \subset C\times \mathbb{R}$, where 
$C=C(\xi,\tau) \equiv 
\{\xi_1 \in \mathbb{R}; |\tau-\xi^3+3\xi\xi_1(\xi-\xi_1)|
\le 2|\tau+\xi^3| \}$.  By using (\ref{eq:Calculus2}) which 
is valid for $b>1/2$,  it is enough to get a
constant upper bound on the following expression 
\begin{equation*}
  \tilde \phi(\xi,\tau) \equiv 
  \frac{|\xi|^2}{\langle \tau+\xi^3\rangle^{-2b'} \langle
   \xi\rangle^{-2s}} \int \frac{|\xi_1(\xi-\xi_1)|^{-2s} 
   \chi_{C(\xi,\tau)}(\xi_1)d\xi_1}
  {\langle \tau-\xi^3+3\xi\xi_1(\xi-\xi_1) \rangle^{2b}}.  
\end{equation*}
From now on we assume $\xi \not =0$.  
By making $\tau=\xi^3(1+y)$, we see that it is sufficient to get an 
upper bound to 
\begin{equation*}
   \tilde \phi_1(\xi,y) \equiv 
  \frac{|\xi|^2}{\langle \xi^3 (y+2)\rangle^{-2b'} \langle
   \xi\rangle^{-2s}} 
  \int \frac{|\xi_1(\xi-\xi_1)|^{-2s} 
   \chi_{C(\xi,\xi^3(1+y))}(\xi_1) d\xi_1}
   {\langle \xi^3 y+3\xi\xi_1(\xi-\xi_1) \rangle^{2b}}.  
\end{equation*}
Now, we make the change of variable $\xi_1=\xi x$.  Since 
$x-x^2=\frac14 - (x-\frac12)^2$, we get
\begin{equation*}
   \tilde\phi_1(\xi,y) \le 
  \frac{|\xi|^{3-4s}}{\langle \xi^3 (y+2)\rangle^{-2b'} \langle
   \xi\rangle^{-2s}} 
  \int \frac{|(x-\frac12)^2-\frac14|^{-2s} \chi_{D_y}(x) dx}
   {\langle \xi^3 (y+3(\frac14-(x-\frac12)^2)) \rangle^{2b}},
\end{equation*}
where $D_y=\{x;|y+3(x-x^2)|\le2|y+2|\}$.  We denote by 
$E_y$ the set given by $\{x;|y+3/4-3x^2|\le 2|y+2|\}$, then we need an 
upper bound on the quantity
\begin{equation*}
  \phi(\xi,y) \equiv 
  \frac{|\xi|^{3-4s}}{\langle \xi^3 (y+2)\rangle^{-2b'} \langle
   \xi\rangle^{-2s}} 
  \int \frac{|x^2-\frac14|^{-2s} 
   \chi_{E_y}(x) dx}
   {\langle \xi^3 (y+3/4-3x^2) \rangle^{2b}}.
\end{equation*}
{\bf{i.)}} First, we suppose $|y+2|>1$.  We remark that  
 $|y-3x^2+\frac34| \le 2|y+2|$ implies $|x^2-\frac14|^{-2s}\le
 c_s|y+2|^{-2s}+ c_s$, for $s \le0$.  If $\phi_1(\xi,y)$ is given by 
(\ref{eq:Estimative-L4}), then we get 
\begin{equation*}
   \phi(\xi,y) \le c_s \phi_1(\xi,y)\Big(1+\frac{1}{|y+2|^{-2s}}\Big)
   \le c_{(s,b,b')}, 
\end{equation*}
where in the last inequality we have used Lemma
\ref{Lemma:Estimative-L4}. \\
{\bf{ii.)}} Now, we assume $|y+2| \le1$.  In $E_y$ we have that 
$|y-3x^2+\frac34| \le 2|y+2|\le 2$, then  
$0\le x^2\le \frac{7}{12}$.  Hence $E_y \subset [-1,1]$.  Moreover, 
$|(y+2)-(3x^2+\frac54)| \le 2|y+2|$ implies 
$|x^2-\frac14|\le \frac54\le 3(x^2+\frac{5}{12})\le 3|y+2|$.    
Therefore,  since $s\le0$, we see that
\begin{equation*}
  \phi(\xi,y) \le 
  \frac{c_s|\xi|^{3-4s}}{\langle \xi^3 (y+2)\rangle^{-2b'} 
  \langle \xi \rangle^{-2s}} 
  \int_0^1 \frac{|y+2|^{-2s} dx}
  {\langle \xi^3 (y+3/4-3x^2) \rangle^{2b}} 
  \le c_{(s,b,b')},
\end{equation*}
where the last inequality is a consequence of (\ref{eq:Estimative-L4}).
\end{proof}

%%%%%%%%%%%%%%%%%%%%%%%%
\begin{Lemma}\label{L5}
If $s \in(-\frac34,-\frac12]$, $b' \in (-\frac12,0]$, and $b>\frac12$
with $b'-b \le \min \{-s-\frac32, s-\frac16\}$, then there
exists $c_{(s,b,b')}>0$ such that
\begin{equation}\label{eq:L5-0}
 \frac{1}{\langle  \tau_1 - \xi_1^3 \rangle^b} \Big( \int\!\!\!\int_B 
 \frac{ |\xi|^{2(1+s)} |\xi \xi_1(\xi-\xi_1)|^{-2s}d\xi d\tau}{\langle
 \xi \rangle^{-2s} \langle \tau + \xi^3 \rangle^{-2b'}  
 \langle \tau-\tau_1-(\xi-\xi_1)^3 \rangle^{2b}}\Big)^{1/2}
 \le c_{(s,b,b')},
\end{equation}
where $B=B(\xi_1,\tau_1)$ is defined as
\begin{equation*}
  B=\{ (\xi,\tau) \in \mathbb{R}^2 ; |\xi_1| \ge 1, |\xi-\xi_1|
  \ge 1, |\tau-\tau_1-(\xi-\xi_1)^3| \le |\tau_1-\xi_1^3|, |\tau+\xi^3| \le
  |\tau_1-\xi_1^3| \}.
\end{equation*}
\end{Lemma}
\begin{proof}We remark that in $B$: 
$|\tau_1+2\xi^3-\xi_1^3-3\xi\xi_1(\xi-\xi_1)| \le 2 |\tau_1-\xi_1^3|$.  
By the inequality (\ref{eq:Calculus2}), it is sufficient to get an
upper bound on the expression
\begin{equation*}
 I(B')= \frac{1}{\langle  \tau_1 - \xi_1^3 \rangle^b} \Big( \int_{B'} 
 \frac{ |\xi|^{2(1+s)} |\xi \xi_1(\xi-\xi_1)|^{-2s}d\xi}{\langle
 \xi \rangle^{-2s}   
 \langle \tau_1+2\xi^3-\xi_1^3-3\xi\xi_1(\xi-\xi_1) \rangle^{-2b'}}\Big)^{1/2},
\end{equation*}
where $B'=\{ \xi \in \mathbb{R}; |\xi_1|\ge 1, |\xi-\xi_1| \ge 1, 
|\tau_1+2\xi^3-\xi_1^3-3\xi\xi_1(\xi-\xi_1)| \le 2|\tau_1-\xi^3| \}$.  
It is not difficult to see that $B'=B_1' \cup B_2'$, where 
$B_1'=\{\xi \in B'; |2\xi^3-3\xi \xi_1(\xi-\xi_1)| 
\le \frac12 |\tau_1 - \xi_1^3| \}$ and 
$B_2'=\{\xi \in B'; \frac12 |\tau_1-\xi_1^3| \le 
|2\xi^3-3\xi \xi_1(\xi-\xi_1)| \le 3 |\tau_1 - \xi_1^3| \}$. \\ 
{\bf{i.)}}  In $B_1'$ we have that:
\begin{equation*}
  \frac12 |\tau_1-\xi_1^3| \le |\tau_1-\xi_1^3 +2\xi^3-3\xi\xi_1(\xi-\xi_1)|,
\end{equation*}
\begin{equation*}
 |\xi| \le |2\xi^3-3\xi\xi_1 (\xi-\xi_1)| \le \frac12|\tau_¹-\xi_1^3|,
\end{equation*}
and
\begin{equation*}
 |\xi \xi_1(\xi-\xi_1)| \le |2\xi^3-3\xi \xi_1(\xi-\xi_1)|
 \le \frac12 |\tau_1-\xi_1^3|.
\end{equation*}
Since $b' \le0$, and $-\frac34 <s \le 0$, it follows that
\begin{eqnarray*}
  I(B_1') &\le&  \frac {c_{b'}}{\langle \tau_1-\xi_1^3  
    \rangle^{b-b'}} \Big( \int_{B_1'} \frac{|\xi|^{2(1+s)}
    |\xi \xi_1(\xi-\xi_1)|^{-2s}d \xi}{\langle \xi \rangle ^{-2s}}  
    \Big)^{1/2}\\
   &\le&  \frac {c_{b'}}{\langle \tau_1-\xi_1^3  
    \rangle^{b-b'+s}} \Big( \int_0^{|\tau_1-\xi_1^3|} 
    (1+\xi)^{2+4s}d \xi \Big)^{1/2}\\
   &\le& \frac{c_{(s,b')}}{\langle \tau_1-\xi_1^3 
    \rangle^{b-b'-\frac32-s}} \le c_{(s,b')},
\end{eqnarray*}
where in the last inequality we have used the fact that 
$b'-b \le -s-\frac32$. \\
{\bf{ii.)}}  First, we remark that in $B_2'$ we have that 
\begin{equation*} 
 3|\xi \xi_1 (\xi -\xi_1)| \le |2\xi^3 -3\xi \xi_1 (\xi - \xi_1)| 
  \le 3 |\tau_1 - \xi_1^3|.  
\end{equation*}
We  define the function 
$\mu(\xi) =\mu_{\xi_1,\tau_1}(\xi)
\equiv \tau_1+2\xi^3-\xi_1^3-3 \xi\xi_1(\xi-\xi_1)$, for
$\xi \in B'$.  We remark that  
$\mu'(\xi) = 3(\xi-\xi_1)^2 +3\xi^2 = 6(\xi-\frac12 \xi_1)^2+\frac32 
\xi_1^2$.  Now, we decompose $B_2'$ into two parts: $B_{2,1}'$ 
and $B_{2,2}'$. \\
Let $B_{2,1}' \equiv \{\xi \in B_2' \,; \, 1 \le |\xi_1| \le 10 |\xi| \}$.
In this set we get:
\begin{equation*}
  1+|\tau_1-\xi_1^3| \le |\xi_1|^3+2|2\xi^3-3\xi \xi_1(\xi-\xi_1)|
  \le c|\xi|^3.
\end{equation*}
Moreover, since $-1 \le s \le -\frac{1}{2}$, it follows from the 
last inequality that
\begin{equation*}
  \frac{|\xi|^{2(1+s)}}{\langle \xi \rangle^{-2s}} \le \frac{1}{\langle \xi
    \rangle^{2(-2s-1)}} 
  \le \frac{c_s}{\langle \tau_1 -\xi_1^3  \rangle^{\frac23 (-2s-1)}}.
\end{equation*}
Since $\mu'(\xi) \ge 3\xi^2$, it follows that 
\begin{equation*}
 \frac{1}{\mu'(\xi)} \le \frac{1}{3\xi^2} \le 
 \frac{c}{\langle \tau_1-\xi_1^3\rangle ^{\frac23}}, \,\,\,\
 \text{for } \xi \in B_{2,1}'.
\end{equation*}
Then, since $b'>-\frac12$, and $-1 \le s \le -1/2$, we get  
\begin{eqnarray*}
  I(B_{2,1}') &\le &  \frac{c_s}{\langle
  \tau_1-\xi_1^3\rangle^{b+\frac13 s}} \Big( \int_{B_{2,1}'} 
  \frac{\mu'(\xi) d \xi}{\langle \mu (\xi)\rangle^{-2b'}} \Big)^{\frac12} \\
  &\le&  \frac{c_s}{\langle
  \tau_1-\xi_1^3\rangle^{b+\frac13 s}} \Big( \int_{|\mu| \le 2|\tau_1-\xi_1^3|} 
  \frac{d \mu}{\langle \mu \rangle^{-2b'}} \Big)^{\frac12}  \\
  &\le& \frac{c_{(s,b')}}{\langle \tau_1-\xi_1^3 \rangle^{b-b'
  +\frac13 s-\frac12}} \le c_{(s,b')},
\end{eqnarray*}
where the last inequality is a consequence of $b'-b \le -s -\frac32$, 
and $s\ge -\frac34$. \\
Finally, we consider $B_{2,2}' \equiv \{ \xi \in  B_2 ' \,; \,
10|\xi| \le |\xi_1|\}$.  Since $-1 \le s \le -\frac12$, we get 
\begin{equation*}
  \frac{|\xi|^{2(1+s)}}{\langle \xi \rangle^{-2s}} \le \frac{1}{\langle \xi
    \rangle^{2(-2s-1)}} \le 1.
\end{equation*}
Moreover, in $B_{2,2}'$ we have that 
\begin{equation*}
  1+|\tau_1-\xi_1^3| \le |\xi_1|^3 +2|2\xi^3-3\xi \xi_1(\xi-\xi_1)| 
  \le c |\xi_1|^3.
\end{equation*}
Since $\mu'(\xi) \ge \frac32 \xi_1^2$, we see that 
\begin{equation*}
 \frac{1}{\mu'(\xi)} \le \frac{c}{\xi_1^2} \le 
 \frac{c}{\langle \tau_1-\xi_1^3\rangle ^{\frac23}}, \,\,\,\
 \text{for } \xi \in B_{2,2}'.
\end{equation*}
Then, by using $b'> -\frac12$, $-1\le s \le -1/2$, and $b'-b \le
s-\frac16$, we see that 
\begin{eqnarray*}
  I(B_{2,2}') &\le &  \frac{c_{b'}}{\langle \tau_1-\xi_1^3 \rangle^{b-b'
  +s-\frac16}} \le c_{b'}.
\end{eqnarray*}
\end{proof}

%%%%%%%%%%%%%%%%%%%%%%%%
\begin{Lemma}\label{L6}
If $b'\le 0$ and $b>\frac12$, then there exists $c_b>0$ such that
\begin{equation}\label{eq:L6-0}
 \frac{|\xi|}{\langle \tau+\xi^3 \rangle^{-b'}} \Big(  \int \!\!\! \int 
 \frac{d\xi_1 d\tau_1}{\langle \tau_1+\xi_1^3 \rangle^{2b} 
  \langle \tau-\tau_1-(\xi-\xi_1)^3 \rangle^{2b}}\Big)^{1/2} \le c_b.
\end{equation}
\end{Lemma}
\begin{proof}
Since $b>\frac12$, it follows from (\ref{eq:Calculus2}) that
\begin{equation*}
  \int \frac{d \tau_1}{\langle \tau_1+\xi_1^3 \rangle^{2b} 
  \langle \tau-\tau_1-(\xi-\xi_1)^3 \rangle ^{2b}} \le
  \frac{c_b}{\langle \tau-\xi^3+3\xi \xi_1(\xi-\xi_1)+2\xi_1^3 
  \rangle^{2b}}.
\end{equation*}
Then, it suffices to prove that 
\begin{equation*}
  \frac{|\xi|^2}{\langle\tau+\xi^3 \rangle^{-2b'}} 
  \int \frac{d \xi_1}{\langle \tau-\xi^3+3\xi \xi_1(\xi-\xi_1) 
  +2\xi_1^3 \rangle^{2b}} \le c.
\end{equation*}
By making the change of variable $\tau=\xi^3(1+z)$, and then 
$\xi_1=\xi x$ inside the integral, it suffices to bound
\begin{equation*}
 \phi(\xi,z)\equiv \frac{|\xi|^3}{(1+|\xi|^3|z+2|)^{-2b'}}
 \int \frac{dx}{(1+|\xi|^3|z+3(x-x^2)+2x^3|)^{2b}}.
\end{equation*}
We define the function $\mu(x)=\mu_z(x)\equiv 2x^3-3x^2+3x+z$.  Then 
$\mu'(x)=6(x-\frac12)^2+\frac32\ge\frac32$.  Since $b'\le0$ and 
$b>\frac12$, it follows that 
\begin{equation*}
 \phi(\xi,z) \le \frac{c|\xi|^3}{\langle \xi^3(z+2)\rangle^{-2b'}}
 \int \frac{\mu_z'(x) dx}{\langle \xi^3 \mu_z (x)\rangle^{2b}}
 \le c |\xi|^3 \int 
 \frac{d\mu}{\langle \xi^3 \mu \rangle^{2b}} =c_b.
\end{equation*}
\end{proof}

%%%%%%%%%%%%%%%%%%%%%%%%
\begin{Lemma}\label{L7}
If $s \in[-\frac34,-\frac12]$, $b' \in [-\frac12,\frac{s}{3}-\frac14]$  
and $b>\frac12$, then there exists $c_{(s,b)}>0$ such that
\begin{equation}\label{eq:L7-0}
 \frac{|\xi|}{\langle  \tau + \xi^3 \rangle^{-b'} \langle \xi
 \rangle^{-s}}  \Big( \int\!\!\!\int_{A_1} 
 \frac{ |\xi_1(\xi-\xi_1)|^{-2s}d\tau_1 d\xi_1}{\langle \tau_1+\xi_1^3
 \rangle^{2b}  \langle \tau-\tau_1-(\xi-\xi_1)^3 \rangle^{2b}}\Big)^{1/2}
 \le c_{(s,b)},
\end{equation}
where $A_1=A_1(\xi,\tau)$ is defined as
\begin{equation*}
  A_1=\{ (\xi_1,\tau_1) \in \mathbb{R}^2 ; |\xi_1| \ge 1, |\xi-\xi_1|
  \ge 1, |\tau-\tau_1-(\xi-\xi_1)^3| \le |\tau_1+\xi_1^3| \le
  |\tau+\xi^3| \}.
\end{equation*}
\end{Lemma}
\begin{proof}
We remark that $A_1 \subset C \times \mathbb{R}$, where 
$C=C(\xi,\tau)= \{ \xi_1 \in \mathbb{R}; 
|\tau-\xi^3+3\xi\xi_1(\xi-\xi_1)+2\xi_1^3| \le 2 |\tau+\xi^3| \}$.  
Since $b>\frac12$, it follows from (\ref{eq:Calculus2}) that it 
is enough to get an upper bound to 
\begin{equation*}
  \frac{|\xi|^2}{\langle \tau+\xi^3\rangle^{-2b'} \langle
   \xi\rangle^{-2s}} \int \frac{|\xi_1(\xi-\xi_1)|^{-2s} 
   \chi_{C(\xi,\tau)}(\xi_1)d\xi_1}
  {\langle \tau-\xi^3+3\xi\xi_1(\xi-\xi_1) +2\xi_1^3 \rangle^{2b}}.  
\end{equation*}
We assume $\xi \not =0$.  Now, we make $\tau=\xi^3(1+y)$, and  
$\xi_1=\xi x$.   Since $s \le0$, it follows that it suffices to bound  
\begin{equation*}
   \phi(\xi,y) \equiv
  \frac{|\xi|^{3-2s}}{\langle \xi^3 (y+2)\rangle^{-2b'}} 
  \int \frac{|x-x^2|^{-2s} \chi_{D_y}(x) dx}
   {\langle \xi^3 (y+3(x-x^2)+2x^3) \rangle^{2b}},
\end{equation*}
where $D_y=\{x;|y+3(x-x^2)+2x^3|\le2|y+2|\}$.  We remark that 
$|x^2-x| \le |2x^3-3x^2+3x-2|$, for all $x \in \mathbb{R}$.  Hence, 
$|x-x^2| \le 3|y+2|$, for $x \in D_y$.  We denote by 
$\mu(x)=\mu_y(x) \equiv 2x^3-3x^2+3x+y$.  \\
{\bf{i.)}} First, we suppose $|y+2|\le1$.  It is not difficult to 
see that $D_y \subset [-1,2]$.  We now assume that $|\xi| \le 1$.  
Since $\mu'(x) \ge \frac32$, $s\le 0$, $b'\le 0$, and $b>\frac12$, 
we have that
\begin{equation*}
  \phi(\xi,y) \le c_s \int_{-1}^2 \frac{|\xi|^3 \mu_y'(x)dx}
  {\langle \xi^3 \mu_y (x) \rangle^{2b}} 
  \le c_s \int \frac{dz}{\langle z \rangle^{2b}} = c_{(s,b)}.
\end{equation*}
Next, we consider the case $|\xi|>1$.   Since $s\le0$ and 
$b'\le 0$, it follows that
\begin{eqnarray*}
 \phi(\xi,y) &\le& c_s \frac{|\xi|^{3-2s}}{|\xi|^{-3b'}|y+2|^{-b'}}
 \int \frac{|y+2|^{-2s}dx}{\langle \xi^3 \mu_y (x)\rangle^{2b}} \\
 &\le& c_s |\xi|^{-2s+3b'} |y+2|^{-2s+b'} 
 \int  \frac{dw}{\langle w \rangle^{2b}} \le c_{(s,b)}, 
\end{eqnarray*}
where the last inequality is a consequence of   
$b>\frac12$, $-\frac12 \le b'\le \frac{s}{3}-\frac14$, 
$-\frac34 \le s \le-\frac14$. \\
{\bf{ii.)}} Finally, we assume $|y+2|>1$.   Since 
$\mu'(x)=6(x^2-x)+3 \ge |x^2-x|$, and $s\le-\frac12$, it follows 
that 
\begin{eqnarray*}
 \phi(\xi,y) 
 &\le& c_s \frac{|\xi|^{3-2s} |y+2|^{-2s-1}}
 {\langle \xi^3 (y+2)\rangle^{-2b'}} 
 \int_{D_y} \frac{ \mu_y'(x)dx}
 {\langle \xi^3 \mu_y (x) \rangle^{2b}} \\
 &\le& c_s (1+|\xi|^3|y+2|)^{-\frac{2s}{3}+2b'} 
 |y+2|^{-\frac{4s}{3}-1} 
 \int_0^{+\infty} \frac{dw}{\langle w\rangle^{2b}}
\end{eqnarray*}
Finally, since $b'\le \frac{s}{3}$, $s \ge -\frac34$, and 
$b>\frac12$, we have that $\phi(\xi,y) \le c_{(s,b)}$.
\end{proof}

%%%%%%%%%%%%%%%%%%%%%%%%
\begin{Lemma}\label{L8}
If $s \in(-\frac34,-\frac12]$, $b' \in (-\frac12,0]$, and $b>\frac12$
with $b'-b \le \min \{-s-\frac32, s-\frac16\}$, then there
exists $c_{(s,b,b')}>0$ such that
\begin{equation}\label{eq:L8-0}
 \frac{1}{\langle  \tau_1 + \xi_1^3 \rangle^b} \Big( \int\!\!\!\int_{B_1} 
 \frac{ |\xi|^{2(1+s)} |\xi \xi_1(\xi-\xi_1)|^{-2s}d\xi d\tau}{\langle
 \xi \rangle^{-2s} \langle \tau + \xi^3 \rangle^{-2b'}  
 \langle \tau-\tau_1-(\xi-\xi_1)^3 \rangle^{2b}}\Big)^{1/2}
 \le c_{(s,b,b')},
\end{equation}
where $B_1=B_1(\xi_1,\tau_1)$ is defined as
\begin{equation*}
  B_1\!=\! \{ (\xi,\tau) \in \mathbb{R}^2 ; |\xi_1| \ge 1, |\xi-\xi_1|
  \ge 1, |\tau-\tau_1-(\xi-\xi_1)^3| \le |\tau_1+\xi_1^3|, |\tau+\xi^3| \le
  |\tau_1+\xi_1^3| \}.
\end{equation*}
\end{Lemma}
\begin{proof}We remark that in $B_1$: 
$|\tau_1+2\xi^3-\xi_1^3-3\xi\xi_1(\xi-\xi_1)| \le 2 |\tau_1+\xi_1^3|$.  
Since $b>\frac12$ and $b'\in [-\frac12,0]$, it follows from 
(\ref{eq:Calculus2}) that it suffices to bound
\begin{equation*}
 I(\tilde B_1)= \frac{1}{\langle  \tau_1 + \xi_1^3 \rangle^b} 
 \Big( \int_{\tilde B_1} 
 \frac{ |\xi|^{2(1+s)} |\xi \xi_1(\xi-\xi_1)|^{-2s}d\xi}{\langle
 \xi \rangle^{-2s}   
 \langle \tau_1+2\xi^3-\xi_1^3-3\xi\xi_1(\xi-\xi_1) \rangle^{-2b'}}\Big)^{1/2},
\end{equation*}
where $\tilde B_1= \tilde B_1(\xi_1, \tau_1)=
\{ \xi \in \mathbb{R}; |\xi_1|\ge 1, |\xi-\xi_1| \ge 1, 
|\tau_1+2\xi^3-\xi_1^3-3\xi\xi_1(\xi-\xi_1)| \le 2|\tau_1+\xi^3| \}$.  
We split $\tilde B_1=\tilde B_{1,1} \cup \tilde B_{1,2}$, where 
\begin{eqnarray*}
  \tilde B_{1,1}&=&\{\xi \in \tilde B_1; |2\xi^3-3\xi \xi_1(\xi-\xi_1)
  -2\xi_1^3| \le \frac12 |\tau_1 + \xi_1^3| \}, \;\;\; \text{and} \\  
  \tilde B_{1,2}&=&\{\xi \in \tilde B_1; \frac12 |\tau_1+\xi_1^3| \le 
  |2\xi^3-3\xi \xi_1(\xi-\xi_1)-2\xi_1^3| \le 3 |\tau_1+\xi_1^3| \}. 
\end{eqnarray*}
{\bf{i.)}}  In $\tilde B_{1,1}$ we have that:
\begin{equation*}
  \frac12 |\tau_1+\xi_1^3| \le |\tau_1-\xi_1^3 +2\xi^3-3\xi\xi_1(\xi-\xi_1)|.
\end{equation*}
Since $2\xi^3-3\xi\xi_1(\xi-\xi_1)-2\xi_1^3 = 
(\xi-\xi_1)(2\xi^2-\xi \xi_1+2\xi_1^2)$,  we also have that
\begin{equation*}
 |\xi|\le |\xi \xi_1(\xi-\xi_1)| \le |2\xi^3-3\xi \xi_1(\xi-\xi_1)-2\xi_1^3|
 \le \frac12 |\tau_1+\xi_1^3|.
\end{equation*}
Since $b' \le0$, and $-\frac34 <s \le 0$, it follows that
\begin{eqnarray*}
  I(\tilde B_{1,1}) 
   &\le&  \frac {c_{b'}}{\langle \tau_1+\xi_1^3  
    \rangle^{b-b'+s}} \Big( \int_0^{|\tau_1+\xi_1^3|} 
    (1+\xi)^{2+4s}d \xi \Big)^{1/2}\\
   &\le& \frac{c_{(s,b')}}{\langle \tau_1+\xi_1^3 
    \rangle^{b-b'-\frac32-s}} \le c_{(s,b')},
\end{eqnarray*}
where in the last inequality we have used the fact that 
$b'-b \le -s-\frac32$. \\
{\bf{ii.)}} In $\tilde B_{1,2}$ we have that 
\begin{equation*} 
 |\xi| \le |\xi \xi_1 (\xi -\xi_1)| \le 
 |2\xi^3 -3\xi \xi_1 (\xi - \xi_1)-2\xi_1^3| 
  \le 3 |\tau_1 + \xi_1^3|.  
\end{equation*}
We  define the function 
$\mu(\xi) =\mu_{\xi_1,\tau_1}(\xi)
\equiv \tau_1+2\xi^3-\xi_1^3-3 \xi\xi_1(\xi-\xi_1)$.  Then  
$\mu'(\xi) = 3(\xi-\xi_1)^2 +3\xi^2 = 6(\xi-\frac12 \xi_1)^2+\frac32 
\xi_1^2$.  We see that $\tilde B_{1,2} = \tilde B_{1,2}^1 
\cup \tilde B_{1,2}^2$, where 
$\tilde B_{1,2}^1 \equiv \{\xi \in \tilde B_{1,2} \,; \, 
1 \le |\xi_1| \le 10 |\xi| \}$, and 
$\tilde B_{1,2}^2 \equiv \{ \xi \in  \tilde B_{1,2} \,; \,
10|\xi| \le |\xi_1|\}$.  The rest of the proof is similar to the proof 
of Lemma \ref{L5}-{\bf{ii.)}}.  Since $-\frac34 \le s \le-\frac12$, 
$b'>-\frac12$ and $b'-b \le -s-\frac32$, it follows that 
$I(\tilde B_{1,2}^1) \le c_{(s,b')}$.  Finally, since 
$-1\le s \le -\frac12$, $b'>-\frac12$ and $b'-b \le s-\frac16$, we 
have that $I(\tilde B_{1,2}^2)\le c_{(s,b')}$. 
\end{proof}

%%%%%%%%%%%%%%%%%%%%%%%%
\begin{Lemma}\label{L9}
If $s \in[-\frac34,-\frac12]$, $b' \in [-\frac12,\frac{s}{3}-\frac14]$  
and $b>\frac12$, then there exists $c_{(s,b)}>0$ such that
\begin{equation}\label{eq:L9-0}
 \frac{|\xi|}{\langle  \tau + \xi^3 \rangle^{-b'} \langle \xi
 \rangle^{-s}}  \Big( \int\!\!\!\int_{A_2} 
 \frac{ |\xi_1(\xi-\xi_1)|^{-2s}d\tau_1 d\xi_1}{\langle \tau_1-\xi_1^3
 \rangle^{2b}  \langle \tau-\tau_1+(\xi-\xi_1)^3 \rangle^{2b}}\Big)^{1/2}
 \le c_{(s,b)},
\end{equation}
where $A_2=A_2(\xi,\tau)$ is defined as
\begin{equation*}
  A_2=\{ (\xi_1,\tau_1) \in \mathbb{R}^2 ; |\xi_1| \ge 1, |\xi-\xi_1|
  \ge 1, |\tau-\tau_1+(\xi-\xi_1)^3| \le |\tau_1-\xi_1^3| \le
  |\tau+\xi^3| \}.
\end{equation*}
\end{Lemma}
\begin{proof}
It is not difficult to see that $A_2\subset C\times \mathbb{R}$, where 
$C=C(\xi,\tau)= \{\xi_1 \in \mathbb{R}; 
|\tau+\xi^3-3\xi_1\xi(\xi-\xi_1)-2\xi_1^3|\le 2 |\tau+\xi^3|\}$.  
Since $b>\frac12$, it follows from (\ref{eq:Calculus2}) that it 
suffices to bound 
\begin{equation*}
  \frac{|\xi|^2}{\langle \tau+\xi^3\rangle^{-2b'} \langle
   \xi\rangle^{-2s}} \int \frac{|\xi_1(\xi-\xi_1)|^{-2s} 
   \chi_{C(\xi,\tau)}(\xi_1)d\xi_1}
  {\langle \tau+\xi^3-3\xi\xi_1(\xi-\xi_1) -2\xi_1^3 \rangle^{2b}}.  
\end{equation*}
We assume $\xi \not =0$.  Then we make $\tau=\xi^3(-1+y)$ and  
$\xi_1=\xi x$.   Since $s \le0$, we see that it is sufficient to bound
\begin{equation*}
   \phi(\xi,y) \equiv
  \frac{|\xi|^{3-2s}}{\langle \xi^3 y\rangle^{-2b'}} 
  \int \frac{|x-x^2|^{-2s} \chi_{E_y}(x) dx}
   {\langle \xi^3 (y-3(x-x^2)-2x^3) \rangle^{2b}},
\end{equation*}
where $E_y=\{x;|y-3(x-x^2)-2x^3|\le2|y|\}$.  We remark that 
$|x^2-x| \le |2x^3-3x^2+3x|$, for all $x \in \mathbb{R}$.  Hence, 
$|x-x^2| \le 3|y|$, for $x \in E_y$.  The rest of the proof is similar 
to the proof of Lemma \ref{L7}.  \\
{\bf{i.)}} If $|y|\le1$, then $E_y \subset [-1,2]$.  First, we suppose 
that $|\xi| \le 1$.  Since $s\le 0$, $b'\le 0$ and $b>\frac12$, it 
follows that $\phi(\xi,y)\le c_{(s,b)}$.  Next, we assume that $|\xi|>1$.   
Since $s \in[-\frac34,-\frac14]$, $b'\in [-\frac12, \frac{s}{3}-\frac14]$ 
and $b>\frac12$, we obtain that $\phi(\xi,y) \le c_{(s,b)}$. \\
{\bf{ii.)}} Since $s \in [-\frac34,-\frac12]$, $b'\le \frac{s}{3}$  and 
$b>\frac12$, we get $\phi(\xi,y) \le c_{(s,b)}$ for $|y|>1$.
\end{proof}

%%%%%%%%%%%%%%%%%%%%%%%%

\begin{Lemma}\label{L10}
If $s \in(-\frac34,-\frac12]$, $b' \in (-\frac12,0]$, and $b>\frac12$
with $b'-b \le \min \{-s-\frac32, \frac{s}{3}-\frac34\}$, then there
exists $c_{(s,b,b')}>0$ such that
\begin{equation}\label{eq:L10-0}
 \frac{1}{\langle  \tau_1 - \xi_1^3 \rangle^b} \Big( \int\!\!\!\int_{B_2} 
 \frac{ |\xi|^{2(1+s)} |\xi \xi_1(\xi-\xi_1)|^{-2s}d\xi d\tau}{\langle
 \xi \rangle^{-2s} \langle \tau + \xi^3 \rangle^{-2b'}  
 \langle \tau-\tau_1+(\xi-\xi_1)^3 \rangle^{2b}}\Big)^{1/2}
 \le c_{(s,b,b')},
\end{equation}
where $B_2=B_2(\xi_1,\tau_1)$ is defined as
\begin{equation*}
  B_2\!=\! \{ (\xi,\tau) \in \mathbb{R}^2 ; |\xi_1| \ge 1, |\xi-\xi_1|
  \ge 1, |\tau-\tau_1+(\xi-\xi_1)^3| \le |\tau_1-\xi_1^3|, |\tau+\xi^3| \le
  |\tau_1-\xi_1^3| \}.
\end{equation*}
\end{Lemma}
\begin{proof}
In $B_2$ we have that  
$| \tau_1+3\xi \xi_1(\xi-\xi_1)+\xi_1^3| \le 2 |\tau_1-\xi_1^3|$.    
Since $b>\frac12$ and $b'\in [-\frac12,0]$, it follows from  
(\ref{eq:Calculus2}) that it is sufficient to bound 
\begin{equation*}
 L(\tilde B_2)= \frac{1}{\langle  \tau_1 - \xi_1^3 \rangle^b} 
 \Big( \int_{\tilde B_2} 
 \frac{ |\xi|^{2(1+s)} |\xi \xi_1(\xi-\xi_1)|^{-2s}d\xi}{\langle
 \xi \rangle^{-2s}   
 \langle \tau_1+3\xi\xi_1(\xi-\xi_1)+\xi_1^3 \rangle^{-2b'}}\Big)^{1/2},
\end{equation*}
where $\tilde B_2=\tilde B_2(\xi_1,\tau_1)
=\{ \xi \in \mathbb{R}; |\xi_1|\ge 1, |\xi-\xi_1| \ge 1, 
|\tau_1+3\xi\xi_1(\xi-\xi_1)+\xi_1^3| \le 2|\tau_1-\xi^3| \}$.  
We see that $\tilde B_2=\tilde B_{2,1} \cup \tilde B_{2,2}$, where 
\begin{eqnarray*}
  \tilde B_{2,1}&=&\{\xi \in \tilde B_2; |2\xi_1^3+3\xi \xi_1(\xi-\xi_1)| 
  \le \frac12 |\tau_1 - \xi_1^3| \}, \;\;\; \text{and} \\  
  \tilde B_{2,2}&=&\{\xi \in \tilde B_2; \frac12 |\tau_1-\xi_1^3| \le 
  |2\xi_1^3+3\xi \xi_1(\xi-\xi_1)| \le 3 |\tau_1-\xi_1^3| \}. 
\end{eqnarray*}
{\bf{i.)}}  In $\tilde B_{2,1}$ we have that
\begin{equation*}
  \frac12 |\tau_1-\xi_1^3| \le |\tau_1+3\xi\xi_1(\xi-\xi_1)+\xi_1^3|, 
  \;\;\; \text{and}
\end{equation*}
\begin{equation*}
 |\xi|\le |\xi \xi_1(\xi-\xi_1)| \le |2\xi_1^3+3\xi \xi_1(\xi-\xi_1)|
 \le \frac12 |\tau_1-\xi_1^3|.
\end{equation*}
Since $b' \le0$, $-\frac34 <s \le 0$ and $b'-b \le -s-\frac32$, it 
follows that $L(\tilde B_{2,1}) \le c_{(s,b')}$. \\
{\bf{ii.)}} In $\tilde B_{2,2}$ we see that 
\begin{equation*} 
 |\xi| \le |\xi \xi_1 (\xi -\xi_1)| \le 
 |2\xi_1^3+3\xi \xi_1 (\xi - \xi_1)| 
  \le 3 |\tau_1 - \xi_1^3|.  
\end{equation*}
We  define the function 
$\mu(\xi) =\mu_{\xi_1,\tau_1}(\xi)
\equiv \tau_1+3 \xi\xi_1(\xi-\xi_1)+\xi_1^3$.  \\
First, we consider $\tilde B_{2,2}^1= \{\xi \in \tilde B_{2,2}; 
\frac{|\xi|}{4} \le|\xi_1| \le 100 |\xi| \}$.  In this set we have that 
\begin{equation*}
 \langle \tau_1-\xi_1^3 \rangle \le |\xi_1|^3 
 +2|2\xi_1^3+3\xi\xi_1(\xi-\xi_1)| \le c|\xi|^3.
\end{equation*}
Since $s \in [-1,-\frac12]$, it follows that 
\begin{equation*}
  \frac{|\xi|^{2(1+s)}}{\langle \xi \rangle^{-2s}} 
  \le \langle \xi \rangle^{2+4s}
  \le c_s \langle \tau_1-\xi_1^3 \rangle^{\frac23 +\frac43 s}.
\end{equation*}
Since $3\xi_1(-\xi_1^3-4\tau_1+4\mu) = (6\xi_1\xi-3\xi_1^2)^2$, 
it follows that  
$|\mu'(\xi)|=|6\xi_1\xi-3\xi_1^2|= \sqrt{3\xi_1(-\xi_1^3-4\tau_1+4\mu)}$.  
Then
\begin{eqnarray*}
 L(\tilde B_{2,2}^1) 
 &\le& \frac{c_s}{\langle \tau_1-\xi_1^3 \rangle^{b+\frac{s}{3}-\frac13}}
 \Big{(} \int_{|\mu| \le 2|\tau_1-\xi_1^3|} 
 \frac{d\mu}{\sqrt{|\xi_1| |-\xi_1^3-4\tau_1+4\mu|} \;\;
 \langle \mu \rangle^{-2b'}} \Big{)}^{\frac12} \\
 &=& \frac{c_s \langle \tau_1-\xi_1^3 \rangle^{-b-\frac{s}{3}+\frac13}}
 {|\xi_1|^{\frac14}} 
 \Big{(} \int_{|\mu| \le 2|\tau_1-\xi_1^3|} 
 \frac{d\mu}{\sqrt{|\xi_1^3/4+\tau_1-\mu|} \;\;
 \langle \mu \rangle^{2(1-(1+b'))}} \Big{)}^{\frac12} \\
 &\le& c_{(s,b')} \frac{ \langle \tau_1-\xi_1^3 
 \rangle^{-b-\frac{s}{3}+b'+\frac56}}
 {|\xi_1|^{\frac14} 
 \langle \xi_1^3+4\tau_1 \rangle^{\frac14}} 
 \le c_{(s,b')} \frac{\langle \tau_1-\xi_1^3 
 \rangle^{-b-\frac{s}{3}+b'+\frac56-\frac{1}{12}}}
 {\langle \xi_1^3+4\tau_1 \rangle^{\frac14}} \\
 &\le& c_{(s,b')} \langle \tau_1-\xi_1^3 
 \rangle^{-b-\frac{s}{3}+b'+\frac34} \le c_{(s,b')},
\end{eqnarray*}
where in the second inequality above we have used (2.11) in \cite{kpv2:kpv2}, and 
$b'>-\frac12$; the last inequality above is a consequence of the 
fact that $b'-b \le \frac{s}{3}-\frac34$. \\
Secondly, we consider 
$\tilde B_{2,2}^2 = \{ \xi \in \tilde B_{2,2}; 1\le |\xi_1| 
\le \frac{|\xi|}{4} \}$.  In this set we have 
\begin{equation*}
 \langle \tau_1-\xi_1^3 \rangle \le c|\xi|^3 
 \;\;\;\;\; \text{and}  \;\;\;\;\;
 \frac{|\xi|^{2(1+s)}}{\langle \xi \rangle^{-2s}} 
  \le c_s \langle \tau_1-\xi_1^3 \rangle^{\frac23 +\frac43 s}.
\end{equation*}
Since $s \in [-1,-\frac12]$ and $b'>-\frac12$, it follows that  
\begin{equation*}
 L(\tilde B_{2,2}^2) 
 \le c_{(s,b')} \frac{ \langle \tau_1-\xi_1^3 
 \rangle^{-b-\frac{s}{3}+b'+\frac56}}
 {|\xi_1|^{\frac14} 
 \langle \xi_1^3+4\tau_1 \rangle^{\frac14}}.
\end{equation*}
If $b'-b-\frac{s}{3}+\frac56 \le 0$, then 
$L(\tilde B_{2,2}^2) \le c_{(s,b')}$.  Thus, we suppose that 
$b'-b-\frac{s}{3}+\frac56 \ge 0$.  Now, we make 
$\tau_1 = \xi_1^3(z-1)/4$.   Hence   
\begin{equation*}
  L(\tilde B_{2,2}^2) \le c_{(s,b')}
  \frac{\langle \xi_1^3(\frac{z-5}{4}) \rangle^{-b-\frac{s}{3}+b'+\frac56}}
  {|\xi_1|^{\frac14} \langle \xi_1^3 z \rangle^{\frac14}}. 
\end{equation*}
Suppose first that $|z|<1$.   Then 
\begin{equation*}
 L(\tilde B_{2,2}^2) \le c_{(s,b,b')}
 \frac{\langle \xi_1^3 \rangle^{b'-b-\frac{s}{3}+\frac56}}
 {|\xi_1|^{\frac14}} \le c_{(s,b,b')} 
 |\xi_1|^{3(b'-b-\frac{s}{3}+\frac56)-\frac14} 
 \le c_{(s,b,b')},
\end{equation*}
where in the last inequality we have used the fact that 
$b'-b \le \frac{s}{3}-\frac34$.  \\ 
Now assume that $|z| \ge 1$.   We see that 
\begin{equation*}
 L(\tilde B_{2,2}^2) 
 \le c_{(s,b,b')} \frac{ \langle \xi_1^3 (z-5)  
 \rangle^{b'-b-\frac{s}{3}+\frac56}}
 {\langle \xi_1^3 z \rangle^{\frac14}}
 \le c_{(s,b,b')}  \langle \xi_1^3 z 
 \rangle^{b'-b-\frac{s}{3}+\frac{7}{12}} 
 \le c_{(s,b,b')},
\end{equation*}
where the last inequality is a consequence of 
$b'-b\le \frac{s}{3}-\frac34$. \\
Finally, we consider the set $\tilde B_{2,2}^3 = \{\xi \in \tilde B_{2,2}; 
100 |\xi| \le |\xi_1| \}$.  In $\tilde B_{2,2}^3$ we have
\begin{equation*}
  \langle \tau_1 - \xi_1^3 \rangle \le |\xi_1|^3 + 
  2|2\xi_1^3+3\xi\xi_1(\xi-\xi_1)| \le c|\xi_1|^3,
\end{equation*}
and $|\mu'(\xi)|=|6\xi_1(\xi-\xi_1)+3\xi_1^2| 
\ge (\frac{3 \times 99}{50}-3) \xi_1^2 \ge 2 \xi_1^2$.  Moreover, 
the fact that $s \in [-1,-\frac12]$ implies that 
\begin{equation*}
  \frac{|\xi|^{2(1+s)}}{\langle \xi \rangle^{-2s}} 
  \le \langle \xi \rangle^{2+4s} \le 1.
\end{equation*}
Then 
\begin{eqnarray*}
  L(\tilde B_{2,2}^3) 
  &\le& \frac{c_s}{\langle \tau_1-\xi_1^3 \rangle^{b+s} \;\; |\xi_1|}
  \Big{(} \int_{\tilde B_{2,2}^3}\frac{|\mu'(\xi)| d\xi}
  {\langle \mu(\xi) \rangle^{-2b'}} \Big{)}^{\frac12}\\ 
  &\le& \frac{c_s}{\langle \tau_1-\xi_1^3 \rangle^{b+s+\frac13}}
  \Big{(} \int_{|\mu|\le 2|\tau_1-\xi_1^3|}\frac{d\mu}
  {\langle \mu \rangle^{-2b'}} \Big{)}^{\frac12} 
  \le c_{(s,b')},
\end{eqnarray*}
where in the last inequality we have used $b'>-\frac12$, 
$b'-b \le \frac{s}{3}-\frac34$ and $s>-\frac34$.
\end{proof}

%%%%%%%%%%%%%%%%%%%%%%%%%%%%%%%%%%%%%%%%%%%%%%%%%
\begin{Remark}\label{Remark:minimum}
It is not difficult to see that 
$\min \{-s-\frac32, \frac{s}{3}-\frac34 \} 
\le \min \{-s-\frac32, s-\frac16 \}$, for $s \ge -\frac78$.
\end{Remark}

The following proposition is the main result of this subsection.
%%%%%%%%%%%%%%%%%%%%%%%%%%%%%%%%%%%%%%%%%%%%%%%%%
\begin{Proposition}\label{prop5}
Given $s>-\frac{3}{4}$,  there exist
$b' \in (-\frac{1}{2},0)$ and $\epsilon_s
>0$ such that for any $b \in(\frac{1}{2}, b'+1]$ with $b'+1-b\le \epsilon_s$
\begin{eqnarray}
 \|(vv)_x\|_{X_{s,b'}^1}\quad &\leq& \quad c_{(s,b,b')}
 \,\|v\|_{X_{s,b}^{-1}}^2, 
 \label{eq:bil1}  \\
 \|(uu)_x\|_{X_{s,b'}^{-1}}\quad &\leq& \quad c_{(s,b,b')}
 \,\|u\|_{X_{s,b}^1}^2,
 \label{eq:bil2} \\
 \|(uv)_x\|_{X_{s,b'}^1}\quad &\leq& \quad c_{(s,b,b')} 
 \,\|u\|_{X_{s,b}^1} \,\|v\|_{X_{s,b}^{-1}}, 
 \label{eq:bil3}  \\
 \|(uv)_x\|_{X_{s,b'}^{-1}}\quad &\leq& \quad c_{(s,b,b')}
 \,\|u\|_{X_{s,b}^1} \,\|v\|_{X_{s,b}^{-1}},
 \label{eq:bil4}  
\end{eqnarray}
where $c_{(s,b,b')}$ is a positive constant depending on $s$, $b$, and $b'$.
\end{Proposition}
\begin{proof} Similar to the proof of Corollary 2.7 in
\cite{kpv2:kpv2}.  Here, we use Lemmas \ref{L3}-\ref{L5} 
to prove (\ref{eq:bil1}) and (\ref{eq:bil2}); the positive number 
$\epsilon_s$ is given by 
\begin{equation*}
 \epsilon_s  =\left\{
  \begin{array}
  [c]{l}%
  \min \{-s-\frac12,s+\frac56 \},\;\;\;\; 
  s\in (-\frac34,-\frac12),\\
  \frac14,\;\;\;\; s\ge0, \\
  \min \{-s'-\frac12,s'+\frac56 \}, \;\;\;\;
  s \in [-\frac12,0),
  \end{array}
  \right.  
\end{equation*}
where $s'$ is any fixed number in the interval 
$(-\frac{3}{4}, -\frac12)$.\\
Lemmas \ref{L6}-\ref{L10} are used to prove (\ref{eq:bil3}) and 
(\ref{eq:bil4}).  Now, we will sketch a proof of (\ref{eq:bil3}).  
We denote by $f$ and $g$ the functions given by 
$f(\xi,\tau):= \langle \tau+\xi^3\rangle^b \langle \xi \rangle^s 
\hat u (\xi, \tau)$, and 
$g(\xi,\tau):= \langle \tau-\xi^3\rangle^b \langle \xi \rangle^s 
\hat v (\xi, \tau)$.  Then 
$\|f\|_{L^2_{\xi} L^2_{\tau}} = \|u\|_{X_{s,b}^1}$, and 
$\|g\|_{L^2_{\xi} L^2_{\tau}} = \|v\|_{X_{s,b}^{-1}}$.  The case 
$s \ge 0$ follows from Lemma \ref{L6} and from the inequality 
$\langle \xi \rangle^s \le \langle \xi_1\rangle^s \langle \xi-\xi_1 \rangle^s$.  
Suppose now that $-\frac34<s< -\frac12$.  Then
\begin{equation*}
\|(uv)_x\|_{X_{s,b'}^1} \le c_{(s,b)} \|u\|_{X_{s,b}^1} \|v\|_{X_{s,b}^{-1}} 
 + c_s \sum_{j=1}^4 \|I_j\|_{L^2_{\xi}L^2_{\tau}},
\end{equation*}
where the first term on the right-hand side of the last inequality corresponds 
to the case when $|\xi_1| \le 1$ or $|\xi- \xi_1| \le1$ (which reduces to the case 
$s=0$), and 
\begin{equation*}
 I_j := \frac{|\xi|}{\langle \tau+\xi^3\rangle^{-b'} \langle \xi \rangle^{-s}} 
 \int\!\!\!\int_{C_j} \frac{|f(\xi_1,\tau_1)| |g(\xi-\xi_1,\tau-\tau_1)|
 |\xi_1(\xi-\xi_1)|^{-s} d\xi_1 d\tau_1}
 {\langle \tau_1+\xi_1^3 \rangle^b 
 \langle (\tau-\tau_1)-(\xi-\xi_1)^3 \rangle^b}, 
\end{equation*}
where 
\begin{eqnarray*}
 C_1&=&\{(\xi_1,\tau_1); |\xi_1|\ge 1, |\xi-\xi_1| \ge1, 
       |(\tau-\tau_1)-(\xi-\xi_1)^3| \le |\tau_1+\xi_1^3| 
       \le |\tau+\xi^3| \}, \\
 C_2&=&\{(\xi_1,\tau_1); |\xi_1|\ge 1, |\xi-\xi_1| \ge1, 
       |(\tau-\tau_1)-(\xi-\xi_1)^3| \le |\tau_1+\xi_1^3|, \\
    & &|\tau+\xi^3| \le |\tau_1+\xi_1³| \}, \\
 C_3&=&\{(\xi_1,\tau_1); |\xi_1|\ge 1, |\xi-\xi_1| \ge1, 
       |\tau_1+\xi_1^3| \le |(\tau-\tau_1)-(\xi-\xi_1)^3|  
       \le |\tau+\xi^3| \}, \\   
 C_4&=&\{(\xi_1,\tau_1); |\xi_1|\ge 1, |\xi-\xi_1| \ge1, 
       |\tau_1+\xi_1^3| \le |(\tau-\tau_1)-(\xi-\xi_1)^3|, \\
    & &|(\tau-\tau_1)-(\xi-\xi_1)^3| \ge |\tau+\xi^3| \}.       
\end{eqnarray*}
The result now follows from Lemmas \ref{L7}-\ref{L10}.  The case  
$s\in [-\frac12,0)$ follows from the last case and from the 
inequality 
$\langle \xi \rangle^{s-s'} |\xi_1(\xi-\xi_1)|^{s'-s} \le c$, which 
holds for $|\xi_1|\ge 1$ and $|\xi-\xi_1|\ge1$, where $s'$ is any  
fixed number belonging to $(-\frac34,-\frac12)$.   Then
\begin{equation*}
 \epsilon_s  =\left\{
  \begin{array}
  [c]{l}
  \min \{-s-\frac12,\frac{s}{3}+\frac14 \},\;\;\;\; 
  s\in (-\frac34,-\frac12),\\
  \frac12,\;\;\;\; s\ge0, \\
  \min \{-s'-\frac12,\frac{s'}{3}+\frac14 \}, \;\;\;\;
  s \in [-\frac12,0).
  \end{array}
  \right.  
\end{equation*}
\end{proof}

%%%%%%%%%%%%%%%%%%%%%%%%%%%%%%%%%%%%%%%%%%%%%%%%%

\begin{Remark}\label{Remark:general-indices, nonequivalence}
{\bf{i.)}} Suppose $a\in \mathbb{R}\setminus \{0\}$.  By making 
similar calculations as in the proof of Proposition \ref{prop4}, 
it follows that Proposition \ref{prop5} still holds if we replace 
the super-indices $1$ by $a$ and $-1$ by $-a$ and the constant 
$c_{(s,b,b')}$ by $c_{(a,s,b,b')}$. \\
{\bf{ii.)}} Consider the bilinear estimate, 
$\|(uv)_x\|_{X_{s,b'}^{-1}} \leq c_{(s,b,b')}
 \|u\|_{X_{s,b}^{-1}} \|v\|_{X_{s,b}^{-1}}$, of Kenig, Ponce, and 
Vega \cite{kpv2:kpv2}.  In the case 
$|\xi_1| \ge 1$ and $|\xi-\xi_1| \ge1$, by symmetry it 
is possible to assume that 
$|\tau -\tau_1-(\xi-\xi_1)^3| \le |\tau_1-\xi_1^3|$ (see 
the proof of Theorem 2.2-\cite{kpv2:kpv2}), and then 
we need to consider only two regions of integration $A$ and $B$ 
(see Lemmas 2.5-\cite{kpv2:kpv2} and 2.6-\cite{kpv2:kpv2} respectively).  
We note, however, that in the proof of (\ref{eq:bil3}) and (\ref{eq:bil4}) 
there is no such symmetry to assume and for this reason the four regions 
of integration $C_1,\ldots, C_4$ (and Lemmas \ref{L7}-\ref{L10})  
were considered.
\end{Remark}

%%%%%%%%%%%%%%%%%%%%%%%%%%%%%%%%%%%%%%%%%%%%%%%%%

%%%%%%%%%%%%%%%%%%%%%%%%%%%%%%%%%%%%%%%%%%%%%%%%%%%%%%%%%%%%%%%%%%%
%Local Well-Posedness to the Gear-Grimshaw System
%%%%%%%%%%%%%%%%%%%%%%%%%%%%%%%%%%%%%%%%%%%%%%%%%%%%%%%%%%%%%%%%%%
\subsection{Local Well-Posedness to the Gear-Grimshaw System}
From now on we consider a cut-off function
$\psi \in C^{\infty}$, such that $0\leq \psi(t) \leq 1$ and
\begin{align*}
 \psi(t)=\left\{ \begin{array}{ll}
 1& \textrm{if \,$|t| \leq 1$},\\
 0& \textrm{if \,$|t| \geq 2$}.
\end{array} \right.
\end{align*}
We define $\psi_{T}(t)\equiv\psi(t/T)$. To prove Theorem
\ref{teox1} we need the following result.

%%%%%%%%%%%%%%%%%%%%%%%%%%%%%
\begin{Proposition}\label{prop3}
Let $s \in \mathbb{R}$, $-\frac12<b' \leq 0\leq b \leq b'+1$,
$T \in [0,1]$, $a \neq 0$.  Then
\begin{align}
  \|\psi_{1}(t)U_a(t)u_{0}\|_{X_{s,b}^a}\quad = & \quad
  c_{(b,\psi)}\,
  \|u_{0}\|_s,  \label{eq1} \\
  \| \psi_{T}(t)\int_{0}^{t}U_a(t-t')F(t',\cdot)dt'\|_{X_{s,b}^a}
  \quad \leq & \quad c_{(b,b',\psi)}\, T^{b'+1-b}
  \|F\|_{X_{s,b'}^a}, \label{eq2}
\end{align}
where
$\widehat{U_a(t)u_0}(\xi)=\exp\{-iat\xi^3\}\hat{u}_0(\xi)$.
\end{Proposition}
\begin{proof}
(\ref{eq1}) is obvious. The proof of (\ref{eq2}) is
practically done in \cite{gtv:gtv}.
\end{proof}
We now prove the following theorem:
\begin{Theorem}\label{teox1}
The IVP (\ref{eq:gg}) with $r=0$  such that
$A=(a_{ij}) \sim a I$ for some $a \not = 0$ is locally well-posed
for data $(u_0,v_0) \in H^s(\mathbb{R}) \times H^s(\mathbb{R})$,
$s>-3/4$.
\end{Theorem}
\begin{proof}
The proof follows from the theory developed by Bourgain
\cite{bou:bou} and Kenig, Ponce and Vega \cite{kpv2:kpv2}. Since
$A \sim aI$, it follows that $a_{11}=a_{22}=a\neq 0$, and
$a_{12}=a_{21}=0$. Let
\begin{align*}
F(u,v)=b_1(uv)_x+b_2uu_x+b_3vv_x, \quad
G(u,v)=b_4(uv)_x+b_5uu_x+b_6vv_x.
\end{align*}
We will consider (\ref{eq:gg}) in its equivalent integral form. Let
$ U_{-a}(t)$ be the unitary group associated with the linear part
of (\ref{eq:gg}). We consider
\begin{equation*}
\Phi(u,v)(t)=(\,\Phi_1(u,v)(t),\,\Phi_2(u,v)(t)\,),
\end{equation*}
where
\begin{align*}
  \Phi_1(u,v)(t)= &\psi(t)\,U_{-a}(t)u_0-\psi_{T}(t)
  \,\int_0^tU_{-a}(t-t')F(u,v)(t')dt',\\
  \Phi_2(u,v)(t)=
  &\psi(t)\,U_{-a}(t)v_0-\psi_{T}(t)\,\int_0^tU_{-a}(t-t')G(u,v)(t')dt'.
\end{align*}
Let $s>-3/4$.  Let $b$, $b'$ be two numbers given by Proposition
\ref{prop4}, such that $\epsilon \equiv b'+1-b>0$.  We will prove that $\Phi(u,v)$
is a contraction in the following space
\begin{align*}
  X_{s,b,a}^M= \{(u,v)\in X_{s,b}^{-a} \times X_{s,b}^{-a}; \,\,
  \|(u,v)\|_{X_{s,b}^{-a} \times X_{s,b}^{-a}}\leq M\},
\end{align*}
where $\|(u,v)\|_{X_{s,b}^{-a} \times
X_{s,b}^{-a}}\equiv \|u\|_{X_{s,b}^{-a}}+\|v\|_{X_{s,b}^{-a}}$.
First we will prove that $\Phi : X_{s,b,a}^M \mapsto
X_{s,b,a}^M$. Let $(u,v)\in X_{s,b,a}^M$.
By using Propositions \ref{prop3}, \ref{prop4} and the
definitions of $F(u,v)$ and $X_{s,b,a}^M$ we get
\begin{align*}
  \|\Phi_1(u,v)\|_{X_{s,b}^{-a}}\quad \leq & \quad C\,\|u_0\|_s+C\,T^\epsilon
  \|F(u,v)\|_{X_{s,b'}^{-a}}\\
  \quad \leq & \quad \frac{M}{4}+ C\,T^\epsilon M^2\leq  \frac{M}{2},
\end{align*}
where we took $M=4C(\|u_0\|_s+\|v_0\|_s)$ and
$CT^\epsilon M=1/4$.  In a similar way we have
\begin{equation*}
  \|\Phi_2(u,v)\|_{X_{s,b}^{-a}}\leq M/2.
\end{equation*}
Therefore $\|\Phi(u,v)\|_{X_{s,b}^{-a} \times X_{s,b}^{-a}}\leq M$. A
similar argument proves that $\Phi$ is a contraction.  We conclude
the proof by a standard argument.
\end{proof}
%%%%%%%%%%%%%%%%%%%%%%%%%%%%%%%%%%%%%%%%%%%%%%%%%

\begin{Remark}\label{Remark:withoutB-type-spaces}
Consider the IVP (\ref{eq:gg}) under the hypothesis of Theorem
\ref{teox1}.  By making the scale change of variables $\tilde u (t,x)
\equiv u(t,a^{1/3}x)$ and $\tilde v (t,x) \equiv v(t,a^{1/3}x)$ we can 
avoid consideration of the modified Bourgain-type spaces $X_{s,b}^a$ 
to prove local-well posedness for data $u_0, v_0 \in
H^s(\mathbb{R})$ for $s>-3/4$.
\end{Remark}

%%%%%%%%%%%%%%%%%%%%%%%%%%%%%%%%%%%%%%%%%%%%%%%%%

\begin{Remark}\label{Remark:No-mixedterms}
Here, we keep the notations of Section  
\ref{subsection:comment}-{\bf{(1)}}.  Suppose that $r=0$ in system 
(\ref{eq:gg}).  Suppose also that 
$A=(a_{ij}) \sim \text{diag}(\alpha_+,\alpha_-)$, 
where $\alpha_+$ and $\alpha_-$ are the eigenvalues of $A$, with 
$\alpha_+, \alpha_- \in \mathbb{R} \setminus \{0\}$, 
$\alpha_+ \not = \alpha_-$. Suppose moreover that the formula 
$\partial_x(u(t)v(t))=\partial_xu(t)v(t)+u(t)\partial_xv(t)$ holds 
for all $t\in[0,T]$ (this is true for example if $s>1/2$, and 
$u(t),v(t) \in H^s(\mathbb{R})$, for all $t \in [0,T]$).   Under 
these assumptions,  we will show that it is possible to obtain 
system (\ref{eq:well}) from system (\ref{eq:gg}), with $C_1(V)V_x$ 
containing only terms of the form $(v_1v_1)_x$, $(v_2v_2)_x$ and 
$(v_1v_2)_x$, where $V=(v_1,v_2)^t$.   
If $a_{12}=a_{21}=0$, there is nothing to prove.  Then, we suppose 
that $a_{12} \not =0$; the case $a_{21}\not = 0$ is similar.  The 
matrices $T$ and $T^{-1}$ are given by
\begin{equation*}
  T=\left( \begin{array}{cc}
  1  & 1 \\
  \frac{\alpha_+-a_{11}}{a_{12}} & \frac{\alpha_--a_{11}}{a_{12}} 
  \end{array} \right), \;\;
  T^{-1}=\frac{a_{12}}{\alpha_+ -\alpha_-}\left( \begin{array}{cc}
  \frac{a_{11}-\alpha_-}{a_{12}}& 1 \\
  \frac{\alpha_+-a_{11}}{a_{12}} & -1 \end{array} \right).
\end{equation*}
Then $V=(\frac{a_{11}-\alpha_-}{\alpha_+ - \alpha_-}u 
+ \frac{a_{12}}{\alpha_+-\alpha_-}v \; , \;
\frac{\alpha_+ -a_{11}}{\alpha_+ - \alpha_-}u 
- \frac{a_{12}}{\alpha_+-\alpha_-}v )^t$.  Now, we see that 
\begin{equation*}
   C_1(V)V_x=\frac{a_{12}}{\alpha_+-\alpha_-}
  \left( \begin{array}{cc}
  av_1+bv_2  & bv_1+cv_2 \\
  dv_1+ev_2  & ev_1+fv_2
  \end{array} \right) \;
  \left( \begin{array}{c}
  \partial_xv_1 \\
  \partial_xv_2
  \end{array} \right),
\end{equation*}
where $a, b, c, d, e, f$ are real constants depending on 
$b_k$, $k=1,\dots, 6$, $a_{i,j}$, $i,j=1,2$, $\alpha_+$ and 
$\alpha_-$.  The result now follows. 
\end{Remark}

%%%%%%%%%%%%%%%%%%%%%%%%%%%%%%%%%%%%%%%%%%%%%%%%%

\begin{Theorem}\label{teox2}
The IVP (\ref{eq:gg}) with $r=0$  such that 
$a_{12}=a_{21}=0$, $a_{11} = -a_{22} \not= 0$ is locally well-posed
for data $(u_0,v_0) \in H^s(\mathbb{R}) \times H^s(\mathbb{R})$,
$s>-3/4$. 
\end{Theorem}
\begin{proof}
Without loss of generality (see Remark 
\ref{Remark:general-indices, nonequivalence}-{\bf{i.)}}), 
we consider the case $a_{11}=-1$ and $a_{22} =1$.   Let
\begin{align*}
F(u,v)=b_1(uv)_x+b_2uu_x+b_3vv_x, \quad
G(u,v)=b_4(uv)_x+b_5uu_x+b_6vv_x.
\end{align*}
We consider
$\Phi(u,v)(t)=(\,\Phi_1(u,v)(t),\,\Phi_2(u,v)(t)\,)$, where
\begin{align*}
  \Phi_1(u,v)(t)= &\psi(t)\,U_1(t)u_0-\psi_{T}(t)
  \,\int_0^tU_1(t-t')F(u,v)(t')dt',\\
  \Phi_2(u,v)(t)=
  &\psi(t)\,U_{-1}(t)v_0-\psi_{T}(t)\,\int_0^tU_{-1}(t-t')G(u,v)(t')dt'.
\end{align*}
Let $s>-3/4$.  Let $b$, $b'$ be two numbers given by Propositions 
\ref{prop4} and \ref{prop5}, with $\epsilon \equiv b'+1-b>0$.
Proceeding in a similar way as in the proof of Theorem
\ref{teox1}, using Propositions \ref{prop4}-\ref{prop3},  it follows that  
$\Phi(u,v)$ is a contraction in the following space
 $$ \mathcal{X}_{s,b}^M = \{(u,v)\in X_{s,b}^1 \times X_{s,b}^{-1}; \,\,
  \|(u,v)\|_{X_{s,b}^1 \times X_{s,b}^{-1}}\leq M\},$$  
$\|(u,v)\|_{X_{s,b}^1 \times
X_{s,b}^{-1}}\equiv \|u\|_{X_{s,b}^1}+\|v\|_{X_{s,b}^{-1}}$, 
$M=4C(\|u_0\|_s + \|v_0\|_s)$ and $CT^{\epsilon}M =\frac14$.
\end{proof}
\noindent
The following result is an immediate consequence of the last theorem.

%%%%%%%%%%%%%%%%%%%%%%%%%%%%%%%%%%%%%%%%%%%%%%%%%
\begin{Corollary}\label{Corollary:teox2}
Let $s>-\frac34$.  Suppose that $r=0$ in (\ref{eq:gg}).   
%be such that the conditions of Remark \ref{Remark:No-mixedterms} 
%are satisfied, then we can get the IVP (\ref{eq:well}) from the IVP  
%(\ref{eq:gg}) with $C_1(V) V_x$ containing only terms of the form 
%$(v_1v_1)_x$, $(v_2v_2)_x$ and $(v_1v_2)_x$, where $V=(v_1,v_2)^t$.  
Suppose also that $A=(a_{ij}) \sim \text{diag}(\alpha_+,\alpha_-)$, 
where $\alpha_+$ and $\alpha_-$ are the eigenvalues of $A$ with 
$\alpha_+, \alpha_- \in \mathbb{R} \setminus \{0\}$, 
$\alpha_+ = - \alpha_-$.  Then the IVP (\ref{eq:gg}) with $r=0$ 
is LWP for data $u_0,v_0 \in H^s(\mathbb{R})$.
\end{Corollary}

%%%%%%%%%%%%%%%%%%%%%%%%%%%%%%%%%%%%%%%%%%%%%%%%%%%%%%%%%%%%%%%%%%
%Future Work
%%%%%%%%%%%%%%%%%%%%%%%%%%%%%%%%%%%%%%%%%%%%%%%%%%%%%%%%%%%%%%%%%%
\subsection{Future Work}
Suppose $a,a' \in \mathbb{R} \setminus \{ 0\}$  
and $|a| \not = |a'|$.  We remark that an interesting problem for 
a future research is to determine whether or not Proposition 
\ref{prop5} is still true when we replace the super-indices $1$ by $a$  
and $-1$ by $a'$.   
%and the constant $c_{(s,b,b')}$ by $c_{(a,a',s,b,b')}$.  
We point out that this result (in general) is not an immediate 
consequence of the calculations we did here for proving Propositions 
\ref{prop4} or \ref{prop5} or from the calculations done in 
\cite{kpv2:kpv2} to prove Corollary 2.7-\cite{kpv2:kpv2}.  This 
result would let us to prove LWP for the Gear-Grimshaw system 
(\ref{eq:gg}) with $r=0$,  when $a_{12}=a_{21}=0$, 
$|a_{11}| \not = |a_{22}|$, and 
$a_{11}, a_{22} \in \mathbb{R} \setminus \{ 0 \}$.   Moreover, if 
this result is true,  
%and under the assumptions of Remark \ref{Remark:No-mixedterms}, 
we also could obtain LWP for  
system (\ref{eq:gg}) with $r=0$, when 
$A=(a_{ij}) \sim \text{diag}(\alpha_+,\alpha_-)$, where 
$\alpha_+$ and $\alpha_-$ are the eigenvalues of $A$ with 
$\alpha_+, \alpha_- \in \mathbb{R} \setminus \{ 0 \}$, 
$|\alpha_+| \not = |\alpha_-|$.

%%%%%%%%%%%%%%%%%%%%%%%%%%%%%%%%%%%%%%%%%%%%%%%%%%%%%%%%%%%%%%%%%%
%Appendix
%%%%%%%%%%%%%%%%%%%%%%%%%%%%%%%%%%%%%%%%%%%%%%%%%%%%%%%%%%%%%%%%%%
\section{Appendix}  Here we prove some properties of  $X_{s,b}^a$-spaces.
%%%%%%%%%%%%%%%%%%%%%%%%%%%%%

\begin{Lemma}\label{L2}
Let $b \geq 0$, $s \in \mathbb{R}$, and $a_0, a_1$ as in Lemma
\ref{L1}. Then for all $a \neq 0$
\begin{equation*}
  X_{s,b}^{a_0} \cap X_{s,b}^{a_1} \subset X_{s,b}^a,   
  \;\;\;\; \text{and}
\end{equation*}
\begin{equation*}
  \|u\|_{X_{s,b}^a} \le c_{(a,a_0,a_1,b)} 
  (\|u\|_{X_{s,b}^{a_0}} + \|u\|_{X_{s,b}^{a_1}}).
\end{equation*}
\end{Lemma}
\begin{proof}[First proof.] 
Let $v$ be an element of $X_{s,b}^{a_0} \cap X_{s,b}^{a_1}$.  Then
\begin{equation*}
  \|v\|_{X_{s,b}^a}^2=\sum_{j=1}^{4}\int_{A_j}\langle \xi \rangle^{2s}\langle
  \tau+a\xi^3 \rangle^{2b}|\hat v(\xi,\tau)|^{2}d\xi d\tau =
  \sum_{j=1}^{4}I_j,
\end{equation*}
where
\begin{align*}
  A_1=\{(\xi,\tau); \xi\ge0, \tau\ge0\},& \quad
  A_2=\{(\xi,\tau); \xi\le0, \tau\le0\},\\
  A_3=\{(\xi,\tau); \xi>0, \tau<0\}, & \quad A_4=\{(\xi,\tau);
  \xi<0, \tau>0\}.
\end{align*}
We consider the case $a>0$, $a_0>0$, and $a_1<0$; a similar
argument works in the other cases.  It is not difficult to prove,
considering regions  $A_1$ and $A_2$,  that
\begin{equation*}
  I_{1}+I_{2} \le 2\Big(1+\frac{a}{a_0}\Big)^{2b}\|v\|_{X_{s,b}^{a_0}}^2.
\end{equation*}
To estimate $I_3$ and $I_4$ we consider
\begin{align*}
  |\tau+a\xi^3|\quad\leq  & \quad
  |\tau+a_0\xi^3|+\Big{|}a_0\xi^3-\frac{a}{|a_1|}\tau\Big{|}
  +a\Big{|}\xi^3-\frac{1}{|a_1|}\tau\Big{|}\\
 \quad \leq & \quad
  |\tau+a_0\xi^3|+\frac{a_0+a}{|a_1|}|\tau+a_1\xi^3|+\frac{a}{|a_1|}|\tau+a_1\xi^3|,
  \end{align*}
therefore
\begin{equation*}
  I_{3}+I_{4} \leq
  c_b\Big{(}1+\frac{a_0+a}{|a_1|}\Big{)}^{2b}(\|v\|_{X_{s,b}^{a_0}}^2
  +\|v\|_{X_{s,b}^{a_1}}^2).
\end{equation*}
{\it{Second proof.}} We claim that for all $x, \tau \in \mathbb{R}$, 
we have 
\begin{equation}\label{eq:estimative-general}
  \frac{1+|\tau+ax|}{(1+|\tau+a_0x|)+(1+|\tau+a_1x|)} 
  \le \left\langle \frac{a-a_0}{a_1-a_0}\right \rangle.
\end{equation}
We will first prove that for all $\xi \ge 0$ and  
$t \in \mathbb{R}$,
\begin{equation}\label{eq:estimativeJ}
 J(\xi,t)\equiv \frac{\xi+|t+a|}{(\xi+|t+a_0|)+(\xi+|t+a_1|)} 
 \le 1+ \frac{| a-a_0|}{|a_1-a_0|}.
\end{equation}
Since 
\begin{equation*}
\xi+|t+a| \le  \xi+|t+a_0+a-a_0| \le \xi+|t+a_0|+|a-a_0|,
\end{equation*}
it follows that 
\begin{equation*} 
  J(\xi,t) \le 1+\frac{|a-a_0|}{(\xi+|t+a_0|)+(\xi+|t+a_1|)} 
           \le 1+\frac{|a-a_0|}{|t+a_0|+|t+a_1|}.
\end{equation*}
Taking $t=x- (a_0+a_1)/2$, $c_0=(a_1-a_0)/2$ and $x=w\,c_0$, we see that 
\begin{equation*}
\frac{1}{|t+a_0|+|t+a_1|} =  \frac{1}{|x-c_0|+|x+c_0|}
 \le  \frac{1}{|c_0|}\frac{1}{|w-1|+|w+1|}.
\end{equation*}
Since the function
\begin{align*}
f(w)=\frac{1}{|w-1|+|w+1|}= \left\{
\begin{array}{ll}
 1/(2w)& \textrm{if \,$w \ge 1$},\\
 1/2 & \textrm{if \,$-1\le w \le 1$},\\
-1/(2w) & \textrm{if \,$w \le -1$}.
\end{array} \right. 
\end{align*}
satisfies $0 \le f(w) \le 1/2$, (\ref{eq:estimativeJ}) follows.  To 
prove  (\ref{eq:estimative-general}) we take 
$\xi=1/|x|$ and $t= \tau/x$ into (\ref{eq:estimativeJ}).  Hence 
\begin{equation*}
  \|u\|_{X_{s,b}^{a}} \le c_b \left \langle \frac{ a-a_0}{a_1-a_0} 
  \right \rangle^b   
  (\|u\|_{X_{s,b}^{a_0}}+\|u\|_{X_{s,b}^{a_1}}).
\end{equation*}
\end{proof}

%%%%%%%%%%%%%%%%%%%%%%%%%%%%%

Thus we can define $X_{s,b}^{a_0,a_1}\equiv X_{s,b}^{a_0} \cap
X_{s,b}^{a_1}$ with norm given by
$\|w\|_{X_{s,b}^{a_0,a_1}}\equiv\|w\|_{X_{s,b}^{a_0}}+\|w\|_{X_{s,b}^{a_1}}$, 
for $b \ge0$, $s \in \mathbb{R}$, and 
$a_0, a_1 \in \mathbb{R} \setminus \{0\}$ such that $a_0 \not = a_1$.
%For $(w_1,w_2) \in X_{s,b}^{a_0,a_1}\times
%X_{s,b}^{a_0,a_1}$ let us define the norm $\|(w_1,w_2)\|_{X_{s,b}^{a_0,a_1}\times
%X_{s,b}^{a_0,a_1}}\equiv\|w_1\|_{X_{s,b}^{a_0,a_1}}+\|w_2\|_{X_{s,b}^{a_0,a_1}}$.

%%%%%%%%%%%%%%%%%%%%%%%%%%%%%
\begin{Corollary}
Let $b \ge0$ and $s \in \mathbb{R}$.
Let $a_0,\ldots,a_3$ be nonzero real numbers such that
$a_0\neq a_1$, $a_2\neq a_3$.  Then
\begin{equation*}
  \mathcal{X}^{s,b} \equiv X_{s,b}^{a_0} \cap X_{s,b}^{a_1}=X_{s,b}^{a_2}
  \cap X_{s,b}^{a_3}.
\end{equation*}
Moreover, there exist constants $c_0 \equiv c_0(a_0, \ldots, a_3, b)$,
$c_1\equiv c_1(a_0, \ldots, a_3, b)>0$, such that
\begin{equation*}
  c_0\,\,\|w\|_{X_{s,b}^{a_0,a_1}}\leq \|w\|_{X_{s,b}^{a_2,a_3}}
  \leq c_1\,\,\|w\|_{X_{s,b}^{a_0,a_1}}.
\end{equation*}
\end{Corollary}
%%%%%%%%%%%%%%%%%%%%%%%%%%%%%

\begin{Remark}
Suppose $b\ge0$ and $s \in \mathbb{R}$.  If $\varphi \in H^b(\mathbb{R})$ 
and $u_0\in H^{s+3b}(\mathbb{R})$, then
$\varphi(t)u_0(x) \in \mathcal{X}^{s,b}$.
\end{Remark}

%%%%%%%%%%%%%%%%%%%%%%%%%%%%%%%%%%%%%%%%%%%%%%%%%%%%%%%%%%%%%%%%%%

\medskip
\noindent
{\bf{Acknowledgements.}}  An earlier version of this paper was written while
the authors had Post-doctoral positions at IMECC-UNICAMP.   The first
author wishes to thank David Lannes and Thierry Colin 
(Universit\'e Bordeaux I and CNRS) for their hospitality and very 
fruitful conversations about Nonlinear Evolution Equations; special thanks 
also go to Jerry Bona (UIC) for very helpful comments.

%%%%%%%%%%%%%%%%%%%%%%%%%%%%%%%%%%%%%%%%%%%%%%%%%
% The Bibliography
%%%%%%%%%%%%%%%%%%%%%%%%%%%%%%%%%%%%%%%%%%%%%%%%%

\medskip

\end{document}